\documentclass[10pt,notitlepage]{article}

\usepackage[utf8]{inputenc}
\usepackage{hyperref}
\usepackage{tikz-cd}
\usepackage{physics}
\definecolor{lightgray}{rgb}{0.63, 0.63, 0.63}
\definecolor{darkgray}{rgb}{0.23, 0.23, 0.23}

\usepackage{charter}
\usepackage{amsmath,amssymb,amsthm}
\usepackage{enumerate}
\usepackage{url}
\usepackage{booktabs}
\usepackage{epigraph}

\usepackage[capitalize]{cleveref}

\usepackage{titlesec}

\titlespacing*{\section}
{0pt}{10.5ex plus 1ex minus .2ex}{5.3ex plus .2ex}
\titlespacing*{\subsection}
{0pt}{5.5ex plus 1ex minus .2ex}{4.3ex plus .2ex}
\titlespacing*{\subsubsection}
{0pt}{3.5ex plus 1ex minus .2ex}{2.3ex plus .2ex}

\titleformat*{\section}{\centering\Large\sffamily}
\titleformat*{\subsection}{\large\sffamily}

\usepackage[left=1.5in, right=1.5in, footskip=1in, bottom=1.5in]{geometry}

\numberwithin{equation}{subsection}

\usepackage{thmtools}

\declaretheoremstyle[
    spacebelow=2\topsep,
    spaceabove=2\topsep,
    headfont=\normalfont\bfseries,
    bodyfont=\itshape,
    postheadspace=\newline,
    qed=${\lrcorner}$,
    headpunct={},
    notebraces={[}{]}
]{breakit}

\declaretheoremstyle[
    spacebelow=2\topsep,
    spaceabove=2\topsep,
    headfont=\normalfont\bfseries,
    bodyfont=\normalfont,
    postheadspace=\newline,
    qed=${\lrcorner}$,
    headpunct={},
    notebraces={[}{]}
]{breakup}

\declaretheorem[numberlike=equation,style=breakit]{theorem}
\declaretheorem[numberlike=equation,style=breakit]{lemma}
\declaretheorem[numberlike=equation,style=breakit]{corollary}

\declaretheorem[numberlike=equation,style=breakup]{definition}
\declaretheorem[numberlike=equation,style=breakup]{example}
\declaretheorem[numberlike=equation,style=breakup]{note}

\DeclareMathOperator{\spec}{Spec}
\DeclareMathOperator{\mspec}{mSpec}
\DeclareMathOperator{\proj}{Proj}

\DeclareMathOperator{\wt}{wt}
\DeclareMathOperator{\van}{\mathbb{V}}
\DeclareMathOperator{\ide}{\mathbb{I}}
\DeclareMathOperator{\rad}{rad}
\DeclareMathOperator{\vanaff}{\mathbb{V}_{\textnormal{aff}}}
\DeclareMathOperator{\ideaff}{\mathbb{I}_{\textnormal{aff}}}

\DeclareMathOperator{\mult}{mult}
\DeclareMathOperator{\GL}{GL}
\DeclareMathOperator{\orbit}{orbit}
\DeclareMathOperator{\im}{im}

\newcommand{\aff}{\mathbb{A}}
\newcommand{\pee}{\mathbb{P}}
\newcommand{\gm}{\mathbb{G}_m}
\newcommand{\ga}{\mathbb{G}_a}
\newcommand{\bb}[1]{\mathbb{#1}}
\newcommand{\noz}{\setminus\{0\}}

\newcommand{\remove}[1]{\widehat{#1}}
\newcommand{\cover}[1]{\overline{#1}}
\newcommand{\othercover}[1]{\widetilde{#1}}

\newcommand{\nn}{\mathbb{N}}
\newcommand{\zz}{\mathbb{Z}}
\newcommand{\cc}{\mathbb{C}}

\newcommand{\kazn}{k_a[x_0,\ldots,x_n]}
\newcommand{\kzny}{k[y_0,\ldots,y_n]}
\newcommand{\kzn}{k[x_0,\ldots,x_n]}
\newcommand{\pazn}{\pee(a_0,\ldots,a_n)}
\newcommand{\padzn}{\pee(a_0',\ldots,a_n')}

\newcommand{\zi}{\{0\}}
\newcommand{\id}{\mathrm{id}}
\newcommand{\prid}{\mathfrak{p}}
\newcommand{\maid}{\mathfrak{m}}
\newcommand*\elide{\textup{[\,\dots]\,}}
\newcommand{\blank}{{-}}

\newcommand{\kathree}{k_a[x_0,x_1,x_2]}
\newcommand{\pathree}{\pee(a_0,a_1,a_2)}
\newcommand{\poly}{S}

\newcommand{\sslash}{/\!\!/}

\begin{document}

    
    \author{Timothy Hosgood}
    \title{An introduction to varieties in weighted projective space}
    \maketitle

    \emph{Written in 2015.}


    \begin{abstract}
        Weighted projective space arises when we consider the usual geometric definition for projective space and allow for non-trivial weights.
        On its own, this extra freedom gives rise to more than enough interesting phenomena, but it is the fact that weighted projective space arises naturally in the context of classical algebraic geometry that can be surprising.
        Using the Riemann-Roch theorem to calculate $\ell(E,nD)$ where $E$ is a non-singular cubic curve inside $\pee^2$ and $D=p\in E$ is a point we obtain a non-negatively graded ring $R(E)=\oplus_{n\geqslant0}\mathcal{L}(E,nD)$.
        This gives rise to an embedding of $E$ inside the weighted projective space $\pee(1,2,3)$.

        To understand a space it is always a good idea to look at the things inside it.
        The main content of this paper is the introduction and explanation of many basic concepts of weighted projective space and its varieties.
        There are already many brilliant texts on the topic (\cite{Reid:2002uy,IanoFletcher:2015wc}, to name but a few) but none of them are aimed at an audience with only an undergraduate's knowledge of mathematics.
        This paper hopes to partially fill this gap whilst maintaining a good balance between `interesting' and `simple'.\footnote{%
            We do not use the word `simple' here in a patronising way, but rather as a way of saying that we will aim to not drift too far into the colossal world of modern abstract algebraic geometry, and to instead merely nod and wave at it as it passes us by.
        }

        The main result of this paper is a reasonably simple degree-genus formula for non-singular `sufficiently general' plane curves, proved using not much more than the Riemann-Hurwitz formula.
    \end{abstract}


    \vspace{-3em}
    \tableofcontents

    \clearpage


\section{Introduction} 
\label{sec:introduction}


    \epigraph{
        My work always tried to unite the true with the beautiful; but when I had to choose one or the other, I usually chose the beautiful.
    }{Herman Weyl}

    \vspace{-2em}
    
    \subsection{Notation, conventions, assumptions, and citations} 
    \label{sub:notation_and_conventions}

    Effort has been made to avoid using too much presupposed knowledge on the behalf of the reader, though a working knowledge of group theory, ring theory, topology, Riemann surfaces, and suchlike is needed, and familiarity with some level of commutative algebra or algebraic geometry, as well as projective geometry, would not go amiss.
    Most `big' theorems or ideas are at least stated before use or mentioned in passing.
    For the sake of not getting too off course, some other sources will be referenced in the main text, usually in lieu of proofs.

    However, as we progress we require more and more tools, and it would be pure folly to try to keep this text entirely self-contained when so many brilliant books and notes have already been written about so many of the subjects which we touch upon in our journey.
    So at the beginning of a section we might state a few results, or a topic, that we assume the reader already has knowledge of, along with a reference for reading up on it if appropriate.
    In general, this paper is aimed at readers of the same level as the author: very early graduate students.

    \bigskip

    All the assumed knowledge of algebraic geometry can be found in \cite{Reid:1988wa} (available online), which is also just a very useful introduction to algebraic geometry as a whole.
    We quite often quote results of affine algebraic geometry (though always try to remember to cite some sort of reference), but usually avoid quoting results of projective algebraic geometry, since this should really be a special case of the things that we're proving here.\footnote{%
        This does change from \cref{sec:curves_in_weighted_projective_space} onwards though, since it is much easier to use some basic facts about projective plane curves than develop everything from scratch.
    }
    Another very enlightening book by the same author is \cite{Reid:1995tu} (also available online) which deals with the commutative algebra side of algebraic geometry.

    There are a few other useful texts to have at hand as a reference (or simply as a better written exposition) for most general algebraic geometry and commutative algebra, as well as the underlying category and scheme theory (if that floats your boat).
    One is \cite{Agrawal:uf}.
    To quote from the \href{https://math.berkeley.edu/~amathew/cr.html}{\texttt{CRing Project website}}:\footnote{\href{https://math.berkeley.edu/~amathew/cr.html}{\url{https://math.berkeley.edu/~amathew/cr.html}}}
    \begin{quotation}
        \emph{The CRing project is an open source textbook on commutative algebra, aiming to comprehensively cover the foundations needed for algebraic geometry at the level of EGA or SGA.
        It is a work in progress.}
    \end{quotation}
    The other is \cite{Vakil:2015wa}, whose purpose is eloquently summarised in the text itself:
    \begin{quotation}
        \emph{
            This book is intended to give a serious and reasonably complete introduction to algebraic geometry, not just for (future) experts in the field.
            \elide
            Our intent is to cover a canon completely and rigorously, with enough examples and calculations to help develop intuition for the machinery. 
            \elide
            We do not live in an ideal world. For this reason, the book is written as a first introduction, but a challenging one.
            This book seeks to do a very few things, but to try to do them well. Our goals and premises are as follows.
        }
    \end{quotation}
    These are not `classic' references in the way that \cite{Hartshorne:1977we} and \cite{Grothendieck:1960tr} are, but they \emph{are} freely available online, which means that even without access to a library they are easy (and free) to get hold of.

    \bigskip

    The main three references, and inspiration, for this paper were \cite{Reid:2002uy,IanoFletcher:2015wc,Dolgachev:1982hz}.

        \subsubsection{List of notation and assumptions} 
        \label{ssub:list_of_notation_and_assumptions}
        
        The following is a list of any notational quirks or assumed conventions used throughout this paper, unless otherwise stated:
        \begin{itemize}
            \item $A\subset B$ means that $A$ is a \emph{proper} subset of $B$, and $A\subseteq B$ means that $A$ is a \emph{not-necessarily proper} subset of $B$;
            \item if $A$ and $B$ are disjoint sets then we (may) write $A\sqcup B$ to mean the union $A\cup B$;
            \item rings are commutative and unital, with $0\neq1$;
            \item $k$ refers to an algebraically closed field;\footnote{%
                We do also sometimes use $k$ as a summation index, or a general integer.
                Hopefully though, it will always be perfectly clear from context which exactly we mean.
            }
            \item $\bb{N}=\{1,2,\ldots\}$, so $0\not\in\bb{N}$;
            \item when we say `graded ring' we mean specifically a $\zz^{\geqslant0}$-graded ring;
            \item for $a,b\in\zz$ we write $a\mid b$ to mean that $a$ divides $b$, and $a\nmid b$ to mean that $a$ does \emph{not} divide $b$;
            \item $\gm=k^{\times}$ is the multiplicative group of a field;
            \item $\ga=k$ is the field considered just as an additive group;
            \item $\mu^n$ is the cyclic group of order $n$, usually realised as $n$-th roots of unity in $k$;
            \item $\omega_n$ is a primitive $n$-th root of unity (so $\mu^n=\langle\omega_n\rangle$);
            \item $R^G$ denotes the subring of $R$ fixed by $G$, where $R$ is a ring and $G$ is a group that acts on $R$;
            \item $k[x_0,\ldots,x_n]$ is the polynomial ring in $n+1$ independent indeterminates, i.e. the $x_i$ are assumed to be algebraically independent, and the same goes for anything similar, such as $k[u,v,w]$, \emph{unless} otherwise explicitly stated;
            \item $f\colon A\twoheadrightarrow B$ means that $f$ is a surjection;
            \item $f\colon A\hookrightarrow B$ means that $f$ is an injection (and so we usually think of it as an embedding of $A$ into $B$);
            \item $(x_1,\ldots,\remove{x_i},\ldots,x_n)=(x_1,\ldots,x_{i-1},x_{i+1},\ldots,x_n)$ -- that is, a hat above an element of a set or ordered n-tuple means ommision of that element;
            \item the phrase `projective space', without either of the words `straight' or `weighted' in front of it, is always taken to refer to both straight \emph{and} weighted projective space\footnote{Obviously both of these will be defined later on, don't worry.}.
        \end{itemize}



    \subsection{The main aims} 
    \label{sub:aims_of_this_paper}

    This paper is meant to serve as a beginner's guide to weighted projective space, its varieties, and some basic algebra that follows on from this.
    In some places the explanations might be quite slow or dense, and the notation seemingly clumsy and complicated, but this is because nearly all of the currently existing literature on these topics is relatively advanced (at least, to an undergraduate or early graduate), so this paper aims to fill a gap in the market, as it were.

    Whenever there has been a decision to make in terms of either leaving some detail out or including it, the latter has almost always been chosen, even if this might not have been the best editorial choice.
    This is because this paper is not intended to be a textbook, where making the reader puzzle out details themselves helps enormously with their eventual understanding of the topic, but instead really as an exercise for the author -- to never leave a proof or explanation as `an exercise for the reader'.
    This was how many of the proofs were originally passed over in the rough notes that came before this paper.
    But after having written the first draft of this paper, I then had to play the role of the reader, which meant having to complete all the `exercises' anyway.

    At the end of the day though, this paper was written by an undergraduate student as a summer project.
    Because of this, the author doesn't have a complete understanding of each and every topic contained within, let alone of the surrounding area of mathematics and the inner workings of the machinery.
    Sometimes questions are left unanswered -- we apologise for this.


    \subsection{An overview of the journey} 
    \label{sub:an_overview_of_the_journey}

    We start off in \cref{sec:weighted_projective_space} with the preliminary definitions of weighted projective space, giving an example or two, and introduce the idea of affine patches to mirror that of straight projective space.

    \Cref{sec:varieties_in_wps} introduces the idea of a \emph{weighted projective variety} and the topology to which it gives rise.
    Again, this is very similar to the usual notion of projective varieties and the Zariski topology on straight projective space, and so it seems natural to consider the idea of some sort of \emph{weighted projective Nullstellensatz}, which we do in \cref{sub:the_weighted_projective_nullstellensatz}.
    Finally in this section we define the coordinate ring of a weighted projective variety, which has the expected definition, bar the \emph{weighting} of the algebraic indeterminates.

    The Nullstellensatz and the coordinate ring let us now do what algebraic geometry always does: study the link between algebra and geometry.
    This is the main motivation of \cref{sec:algebra}.
    We use the Nullstellensatz to explain the $\proj$ construction and gain some intuition with this algebraic definition of weighted projective space, noting the similarities to that of the $\spec$ construction in affine algebraic geometry.
    Using the idea of \emph{truncation of a graded ring} we show that we can often reduce a weighted projective space to a simpler one, where any $(n-1)$ of the weights are coprime.
    We give these `fully simplified' weighted projective spaces the nice name of `well-formed'.
    To ensure that we don't get lost in the heady realms of abstract algebra, this section ends with a few worked examples and applications.

    Sticking with the theme of following a standard algebraic geometry course we look next, in \cref{sec:curves_in_weighted_projective_space}, at \emph{plane curves}, which are varieties in weighted projective 2-space given by the vanishing of a single weighted-homogeneous polynomial: $X=\van(f)\subset\pee(a_0,a_1,a_2)$.
    We show that sufficiently nice plane curves are also Riemann surfaces, and then use some of the associated machinery (such as the Riemann-Hurwitz formula) to see what else we can find out.
    This section culminates with a version of the degree-genus formula for weighted projective plane curves.

    Finally, \cref{sec:the_view_from_where_we_ve_got_to} is a brief introduction to a few further ideas that follow on from the theory that we have developed up until this point.
    We cover (speedily and not overly rigorously) some ideas about the Hilbert polynomial of a variety and the Hilbert Syzygy Theorem.
    It is probably the most interesting of all sections, but is purely expository, with references to various sources that deal with the subject matter in more detail and depth.
    

    \subsection{Acknowledgements} 
    \label{sub:acknowledgements}

    Many thanks go to Balázs Szendrői for suggesting the subject of, and then supervising, this project.
    His patience in answering my multitude of questions and explaining so many things was the main reason that the author was able to write this paper.

    Thanks must also be given to Edmund de Unger Academic Purposes fund from Hertford College, Oxford, without which the author would not have been able to live in Oxford throughout the composition of this paper.


    \subsection{Corrections} 
    \label{sub:corrections}
    
    It is absurdly unlikely that there are no mistakes in this paper.
    If you find any, or have any recommendations for rewording or reordering of sections, please do email the author at \href{mailto:timhosgood@gmail.com?subject=RE:%20An%20Introduction%20to%20Weighted%20Projective%20Space}{\texttt{timhosgood@gmail.com}} including the version date (as found on the title page) in your message.
    Source files for this document are at \href{https://github.com/thosgood/introduction-to-wps}{\texttt{github.com/thosgood/introduction-to-wps}} and are available for use and modification.


    \bigskip

    Apologies are made in advance for the often overly-florid language that is most likely not-at-all suitable for rigorous mathematics.


\section{Weighted projective space} 
\label{sec:weighted_projective_space}

At a first glance, it's very possible that the idea of \emph{weighted projective space} seems like a needless generalization of the usual projective space (herein referred to as \emph{straight projective space}) to which we have all grown to know and love, especially when we show that we can simply embed weighted projective space into a large enough straight projective space.
But it turns out that, quite often, using a weighted projective space instead of its embedding results in much simpler equations and a generally more manageable beast.
So some projective varieties can be more easily described as weighted projective spaces (similar to how some seemingly complicated varieties can turn out to be just a Veronese embedding of some smaller straight projective space), but there is also the fact that weighted projective space itself encompasses the idea of straight projective space and allows us to study interesting ideas in more generality.
We cover some particularly natural uses and applications of weighted projective space in \cref{sub:elliptic_curves_and_friends}.

Like in a lot of algebraic geometry, there are two main ways of approaching a topic: geometrically and algebraically.
The geometric way below mirrors the method usually used in approaching projective space and is a very `hands-on' construction.
The algebraic way uses some language from very basic scheme theory, dealing with the $\proj$ construction, and is covered in \cref{sec:algebra}.

    \subsection{Geometric construction of weighted projective space} 
    \label{sub:geometric_construction}

    \begin{definition}[Weights and their induced action]
        Let $a=(a_0,\ldots,a_n)$ with $a_i\in\bb{N}$, and define the corresponding action of $\gm$ (which we write as $\gm^{(a)}$ to avoid confusion) on $\aff^{n+1}\noz$ by
        \begin{equation}
            \lambda\cdot(x_0,\ldots,x_n) = (\lambda^{a_0}x_0,\ldots,\lambda^{a_n}x_n).
        \end{equation}
        We call $a_1,\ldots,a_n$ the \emph{weights}.    
    \end{definition}

    We use this action to define weighted projective space as follows.

    \begin{definition}[Weighted projective space]
    \label{def:geometric_wps}
        Let $a=(a_0,\ldots,a_n)$ be a weight.
        Define \emph{$a$-weighted projective space} as the quotient
        \[
            \pee(a_0,\ldots,a_n) = (\aff^{n+1}\noz)/\gm^{(a)}.
        \]
        
        We write points in $\pee(a_0,\ldots,a_n)$ as $|x_0:\ldots:x_n|_a$, which represents the equivalence class of the point $(x_0,\ldots,x_n)\in\aff^{n+1}\noz$, omitting the subscript $a$ if it is clear that we are working in $\pee(a)=\pee(a_0,\ldots,a_n)$.
    \end{definition}

    It is easier to deal with weighted projective space once we have some technical lemmas and another perspective under our belt, and so the majority of big examples will come at the end of a later section, but here are a few right now to help build some kind of intuition.

    \begin{example}[Recovering straight projective space]
        If we set all $a_i$ to be $1$ then the above definition coincides with that of straight projective space:
        \[
            \pee(1,\ldots,1) = (\aff^{n+1}\noz)/\gm = \pee^n
        \]
        where $|x_0:\ldots:x_n|\sim\lambda|x_0:\ldots:x_n|=|\lambda x_0:\ldots:\lambda x_n|$, and we more commonly write $[x_0:\ldots:x_n]$ for the coordinates.
    \end{example}

    \begin{example}[Something that looks a bit like a cone]
    \label{ex:something_like_a_cone}
        \textit{We cover this example in more detail later, since it turns out to be an interesting one, so here we look at it quite briefly and not particularly rigorously.}
        
        Consider $\pee(1,1,2)$.
        Just like in straight projective space, any point is invariant under scaling, but here with respect to the weighting.
        For example,
        \[
            |1:0:2| = 3\cdot|1:0:2| = |3:0:18|.
        \]

        Now we try to get a bit more of a grasp on the space as a whole.
        Define a map
        \[
            \varphi\colon |x_0:x_1:x_2| \mapsto [x_0^2:x_0x_1:x_1^2:x_2].
        \]
        We claim that this map has its image in $\pee^3$.
        Firstly, since at least one of the $x_i$ is non-zero (by the definition of weighted projective space), at least one of the monomials will also be non-zero.
        But then all that we need to check is that the image is invariant under scaling, i.e. that
        \[
            \lambda\cdot[x_0^2:x_0x_1:x_1^2:x_2] = [\lambda x_0^2:\lambda x_0x_1:\lambda x_1^2:\lambda x_2] = [x_0^2:x_0x_1:x_1^2:x_2].
        \]
        Simply using definitions we get that
        \begin{align*}
            [x_0^2:x_0x_1:x_1^2:x_2] &= \varphi(|x_0:x_1:x_2|) \\
            &= \varphi\left(\lambda^{\frac12}\cdot|x_0:x_1:x_2|\right) \\
            &= \varphi\left(|\lambda^{\frac12}x_0:\lambda^{\frac12}x_1:\lambda x_2|\right) = [\lambda x_0^2:\lambda x_0x_1:\lambda x_1^2:\lambda x_2].
        \end{align*}

        Even though we haven't really defined what an isomorphism should be for weighted projective spaces, it makes sense to think that, if we can find an inverse map that is also given by polynomials in each coordinate, then we can think of $\pee(1,1,2)$ and its image under $\varphi$ in $\pee^3$ as isomorphic.
        That is, we can think of $\varphi$ as an embedding of $\pee(1,1,2)$ in $\pee^3$.

        To construct our inverse map we take some point $[y_0,y_1,y_2,y_3]$ in the image.
        Sadly, even though $k$ is algebrically closed, we can't take $|y_0^{1/2}:y_2^{1/2}:y_3|$ as our inverse map, since this is not polynomial in each coordinate.
        But we do know that
        \[
            y_0 = x_0^2,~~ y_1 = x_0x_1,~~ y_2 = x_1^2,~~ y_3 = x_2
        \]
        for some $|x_0:x_1:x_2|\in\pee(1,1,2)$, and so
        \begin{align*}
            |x_0:x_1:x_2| =
            \begin{cases}
                x_0\cdot|x_0:x_1:x_2| = |x_0^2:x_0x_1:x_0^2x_2| =|y_0:y_1:y_0y_3|; \\
                x_1\cdot|x_0:x_1:x_2| = |x_0x_1:x_1^2:x_1^2x_2| =|y_1:y_2:y_2y_3|,
            \end{cases}
        \end{align*}
        where we choose whichever option gives us a point in $\pee(1,1,2)$, i.e. depending on whether or not all of $y_0,y_1,y_3$ are zero.

        So we have two mutually inverse polynomial maps, which we (for the moment) are content with calling an isomorphism:
        \begin{equation*}
            \begin{array}{rlcl}
                \varphi\colon & \pee(1,1,2) & \to & X\subset\pee^3 \\
                & |x_0:x_1:x_2| & \mapsto & [x_0^2:x_0x_1:x_1^2:x_2] \\
                & & & \\
                \varphi^{-1}\colon & X &\to & \pee(1,1,2) \\
                & |y_0:y_1:y_2:y_3| & \mapsto &
                    \begin{cases}
                        |y_0:y_1:y_0y_3|\quad\mbox{if}~y_0,y_1,y_3\neq0; \\
                        |y_1:y_2:y_2y_3|\quad\mbox{otherwise}.
                    \end{cases}
            \end{array}
        \end{equation*}
        Then understanding $\pee(1,1,2)$ becomes a matter understanding the set $X\subset\pee^3$.
        We don't know yet what properties $X$ has exactly, but we will find out later on.

        One thing to notice though is that, on the patch $\{x_2=0\}$ of $X\subset\pee^3$, we have something that looks a lot like the Veronese embedding of $\pee^1$ into $\pee^2$: the \emph{rational normal curve of degree $2$}, also known as the \emph{flat conic}.\qedhere
    \end{example}

    The way that we came up with the idea of the map in \cref{ex:something_like_a_cone} was to construct all the degree two monomials, but where we consider $x_2$ as being a degree two element already.
    Generally, though, we see that it looks like we should be able to embed weighted projective space into some straight projective space of high enough dimension, using something that looks remarkably like a Veronese embedding.

    It's not entirely clear why exactly this should work at this point, but it turns out to be a much easier idea to justify once we have introduced the notion of weighted projective spaceusing the $\proj$ construction, and truncation, later on.


    \subsection{Coordinate patches} 
    \label{sub:coordinate_patches}

    On straight projective space, picking some $0\leqslant i\leqslant n$, we have the standard decomposition
    \[
        \pee^n = \{[x_0:\ldots:x_n]\mid x_i=0\}\sqcup\{[x_0:\ldots x_n]\mid x_i\neq0\} \cong \pee^{n-1}\sqcup\aff^n
    \]
    where the sets $H_i=\{x_i=0\}$ are closed, and so the $U_i=\{x_i\neq0\}$ are open.%
    \footnote{%
        Here we mean in the Zariski topology, but it is also true in the quotient topology coming from the definition as \mbox{$\pee^n=(\aff^n\setminus\{0\})/\gm$}.%
    }

    It makes sense to consider the same sets $H_i$ and $U_i$ in weighted projective space, and even though we have yet to really define a Zariski-style topology%
    \footnote{%
        This does have the expected definition: simply define closed sets to be those that can be written as $\van(I)$ for some ideal $I$ in a suitable polynomial ring.
        But we do have to worry about which polynomial ring, and what it means to evaluate a polynomial at a point in wps.
        This is all covered in \cref{sub:weighted_projective_varieties}, in particular in \cref{defn:zariski-topology}.%
    } it turns out that these sets are closed and open (respectively) as before, in both the Zariski and quotient topologies.

    In this decomposition the $U_i$, often called \emph{coordinate patches} or \emph{affine patches}, turn out to be very useful, since they cover the whole of projective space and are isomorphic to affine space, which is much easier to visualise geometrically in most cases.
    So given some projective variety, we can see how it intersects with the affine patches $U_i$ and study these using all our familiarity with affine space.

    But in weighted projective space we have a slight issue, namely that the $U_i$ are \emph{not} isomorphic to $\aff^n$, but instead some quotient of $\aff^n$ by a finite group.
    They still deserve the name `affine patches' then, since the quotient of an affine variety by a finite group is again an affine variety%
    \footnote{%
        See \cite[Definition~3.6,~Theorem~3.8]{Hoskins:2012uq}.
    } but a little bit more work needs to be done to see what exactly they look like.

    \begin{definition}[Quotient of affine space by a cyclic group]\label{defn:quotient-affine-cyclic}
        Define an action of $\mu^{a_i}$ on $\aff^n$, called the \emph{action of type $\frac{1}{a_i}(a_0,\ldots,\remove{a_i},\ldots,a_n)$}, by
        \begin{equation}
            \omega_{a_i}\cdot(x_0,\ldots,\remove{x_i},\ldots,x_n) = (\omega_{a_i}^{a_0}x_0,\ldots,\remove{\omega_{a_i}^{a_i}x_i},\ldots,\omega_{a_i}^{a_n}x_n).
        \end{equation}
        This induces an action of $\mu^{a_i}$ on $k[x_0,\ldots,\remove{x_i},\ldots,x_n]$ given by $\omega_{a_i}\cdot x_j = \omega_{a_i}^{a_j}x_j$ and thus\footnote{%
            Again, see \cite[Definition~3.6,~Theorem~3.8]{Hoskins:2012uq}.
        } gives rise to the affine quotient variety
        \[
            \aff^n/\mu^{a_i} = \mspec \left(k[x_0,\ldots,\remove{x_i},\ldots,x_n]^{\mu^{a_i}}\right)
        \]
        as well as the quotient map $\pi_i = (\iota_i)_{\#}\colon\aff^n\to\aff^n/\mu^{a_i}$ corresponding to the inclusion
        \[
            \iota_i\colon k[x_0,\ldots,\remove{x_i},\ldots,x_n]^{\mu^{a_i}}\hookrightarrow k[x_0,\ldots,\remove{x_i},\ldots,x_n].\qedhere
        \]
    \end{definition}

    \begin{lemma}[Affine patches in wps]\label{lem:affine-patches}
        With $U_i=\{|x_0:\ldots:x_n|\in\pee(a_0,\ldots,a_n) : x_i\neq0\}$ and the quotient $\aff^n/\mu^{a_i}$ defined as in \cref{defn:quotient-affine-cyclic} we have
        \[
            U_i \cong \aff^n/\mu^{a_i}
        \]
        where we mean isomorphic in the usual sense: there exists an algebraic morphism given by a polynomial map with polynomial inverse.
        We often write $A_i=\aff^n/\mu^{a_i}$.
    \end{lemma}

    \begin{proof}
        We can write every point in $U_i$ as $|x_0:\ldots:1:\ldots:x_n|$ where the $1$ is in the $i$-th place.
        Write a point in $A_i$ as $[(y_1,\ldots,y_n)]$, that is, the equivalence class of $y=(y_1,\ldots,y_n)\in\aff^n$ consisting of the points in the $\mu^{a_i}$-orbit of $y$.
        Define the morphism $\varphi_i\colon U_i\to A_i$ by
        \[
            \varphi_i\colon |x_0:\ldots:1:\ldots:x_n| \mapsto \big[(x_0,\ldots,\remove{1},\ldots,x_n)\big]
        \]
        which we need to show is a well-defined map.
        Our concern is that this map might rely on our choice of $x_i$.
        By definition, where $\omega_i=\omega_{a_i}$ is a primitive $a_i$-th root of unity,
        \[
            |x_0:\ldots:1:\ldots:x_n| = |\omega_i^{a_0}x_0:\ldots:\omega_i^{a_i}:\ldots:\omega_i^{a_n}x_n| = |\omega_i^{a_0}x_0:\ldots:1:\ldots:\omega_i^{a_n}x_n|.
        \]
        Now
        \begin{align*}
                \varphi_i\big(|\omega_i^{a_0}x_0:\ldots:1:\ldots:\omega_i^{a_n}x_n|\big) &= \big[(\omega_i^{a_0}x_0,\ldots,\remove{1},\ldots,\omega_i^{a_n}x_n)\big] \\
                &= \big[\omega_i\cdot(x_0,\ldots,\remove{1},\ldots,x_n)\big] = \varphi_i\big(|x_0:\ldots:1:\ldots:x_n|\big),
        \end{align*}
        and so our map is well defined.
        Define its inverse morphism $\varphi_i^{-1}\colon A_i\to U_i$ by
        \[
            \varphi_i^{-1}\colon \big[(y_1,\ldots,y_n)\big] \mapsto |y_1:\ldots:y_{i-1}:1:y_i:\ldots:y_n|
        \]
        and note that it is indeed an inverse to $\varphi$.
        We can show that it is well defined in exactly the same way that we did for $\varphi_i$, and clearly both $\varphi_i$ and $\varphi_i^{-1}$ are polynomial in each coordinate.
    \end{proof}

    Although the $U_i$ give us a nice affine space to work with, it is a quotient space and so can be quite tricky to realise at times.
    Much easier is the idea of looking at the covering space of the affine patches, since it turns out that if the covering space of the affine patch has certain nice properties then so too does the affine patch.
    We now define a few terms to avoid confusion after all this talk of `affine patches'.

    \begin{definition}[Quotient and covering affine patches]\label{defn:quotient-covering-affine-patches}
        Given some subset $X\subseteq\pee(a_0,\ldots,a_n)$ we define the \emph{quotient affine patches $X_i\subseteq A_i$} and the \emph{covering affine patches $\cover{X_i}\subseteq\aff^n$} of $X$ as
        \[
            \begin{array}{rllll}
                X_i &=& X\cap U_i &\subseteq& A_i=\aff^n/\mu^{a_i} \\[0.45em]
                \cover{X_i} &=& \pi_i^{-1}(X\cap U_i) &\subseteq& \aff^n
            \end{array}
        \]
        where we use the isomorphism \mbox{$U_i\cong A_i$} and the quotient map $\pi_i\colon \aff^n\to A_i$.
    \end{definition}

    So we have two different types of affine patches to think of: those that we glue together to get the ambient weighted projective space, which are in some way `folded up' (the quotient affine patches); and those that come from `unfolding' the aforementioned ones (the covering affine patches).
    The reason for considering both is that they can be equally useful, but in different ways.
    Really, the covering affine patches come into their own in \cref{sec:curves_in_weighted_projective_space}, since we already know many useful facts about varieties in $\aff^n$.



\section{Weighted projective varieties} 
\label{sec:varieties_in_wps}

\emph{Although we are dealing with so-called weighted-homogeneous ideals here, we are really just looking at homogeneous ideals (i.e. graded submodules where we consider the ring as a module over itself) of a graded ring.
Thus most standard proofs can be used for the vast majority of the lemmas in this section, and at times we simply refer to them instead of providing our own.
A good reference is \cite[Chapter~6,~Section~1]{Agrawal:uf}.}

\bigskip

Before we define the notion of \emph{weighted projective varieties} we need to cover a few formalities.
The main one is `how do we define evaluating a polynomial at a point in weighted projective space, and is this well defined?'.
It turns out that, just as in straight projective space, evaluating a polynomial at a point is \emph{not} well defined, but seeing whether or not a point is a zero of a polynomial \emph{is}, as long as our polynomial is homogeneous, but in a slightly different sense.
We start off this section in quite a dull manner, with quite a few definitions in a row, and most of them quite expected or natural, but then play around with them to see what we can get.

\begin{definition}[Weighted polynomial ring]
    Define the \emph{polynomial ring in $n+1$ variables with weighting $a=(a_0,\ldots,a_n)$} as
    \[\kazn\text{ with }\wt x_i=a_i.\]
    That is, we think of $x_i$ as a degree $a_i$ monomial and thus, for example,
    \[\deg\left(\prod_{i=0}^n x_i^{c_i}\right) = \sum_{i=0}^n a_i c_i.\]
    \emph{If we omit the subscript $a$ and simply write $k[x_0,\ldots,x_n]$ then we mean the polynomial ring in $n+1$ variables with the usual grading, i.e. $a_i=1$ for all $i$.}

    We sometimes use the phrase \emph{weighted degree} to be clear that we are including the weighting in our calculation of the degree.
    Note that $\deg\lambda=0$ for any $\lambda\in k$.
    For a general polynomial $f\in\kazn$ we define the \emph{degree} $\deg f$ as the maximum of all the degrees of the monomials in $f$.\footnote{%
        Just as we do for polynomials normally.
        So, for example, $f=x^2+y^2\in k_{(1,2)}[x,y]$ has $\deg f=4$.
    }
\end{definition}

\begin{note}
    An important thing to note is that this weighting changes the \emph{grading} of the ring, but it doesn't change the underlying $k$-algebra structure.
    \emph{So $\kazn$ is Noetherian.}
    Further, it means that until \cref{sec:algebra}, when we start talking about the graded ring structure, our choice of notation $\kazn$ vs. $k[x_0,\ldots,x_n]$ is reasonably arbitrary, and doesn't particularly matter.
\end{note}

\begin{definition}[Weighted-homogeneous polynomial]\label{defn:wh-poly}
    Let $f\in k[x_0,\ldots,x_n]$ where $\wt x_i=a_i$ for some weight $a=(a_0,\ldots,a_n)$.
    We say that $f$ is \emph{$a$-weighted-homogeneous of degree $d$}\footnote{%
        Or simply \emph{weighted-homogeneous of degree $d$} if it is clear that we are working with the weight $a$.
    } if each monomial in $f$ is of weighted degree $d$, i.e. there exist $c_i\in k$ and some $m\in\nn$ such that
    \[f = \sum_{i=1}^m c_i\left(\prod_{j=0}^n x_j^{d_j^{(i)}}\right)\]
    and, for all $0\leqslant i\leqslant n$,
    \[\sum_{j=0}^n a_j d_j^{(i)} = d.\]

    We write $\kazn_d\subset\kazn$ to mean the additive group of all weighted-homogeneous polynomials of degree $d$.
\end{definition}

Note that if $f$ is $a$-weighted-homogeneous of degree $d$ then for any $\lambda\in\gm$ then, by \cref{defn:wh-poly},
\begin{equation}\label{eq:wh-poly}
    f(\lambda^{a_0}x_0,\ldots,\lambda^{a_n}x_n) = \lambda^d f(x_0,\ldots,x_n).
\end{equation}

So let $p=|p_0:\ldots:p_n|\in\pazn$ and $f\in \kazn$.
By definition we also have that $p=|\lambda^{a_0} p_0:\ldots:\lambda^{a_n} p_n|$ for any $\lambda\in\gm$.
In particular we can assume that $\lambda\neq1$.
But then, using \cref{eq:wh-poly},
\[f(\lambda^{a_0}p_0,\ldots,\lambda^{a_n}p_n) = f(p_0,\ldots,p_n) \quad\text{if and only if}\quad f(p_0,\ldots,p_n)=0.\]
Thus the idea of evaluating an $a$-weighted-homogeneous polynomial $f$ at a point $p\in\pee(a)$ doesn't make sense in general, but looking at the points $p\in\pee(a)$ at which $f$ vanishes \emph{does} make sense.
That is, it is well defined to write that $f(p)=0$ for some $f\in \kazn$ and $p\in\pee(a)$.
We will come back to this point shortly, in \cref{sub:weighted_projective_varieties}.

\begin{definition}[Weighted-homogeneous ideal]
    We say that an ideal $I\triangleleft \kazn$ is \emph{$a$-weighted-homogeneous}\footnote{%
        From now on we will stop pointing out that we might sometimes omit the weight $a$ if it is clear which weight we are working with, but this is still the case.
        The convention that the author has tried to stick to is to simply say `weighted-homogeneous' if we have already said that $I\triangleleft\kazn$, since then the weight $a$ is clear.%
    } if it is generated by $a$-weighted-homogeneous elements (of not necessarily the same degree).
\end{definition}

\begin{example}
    Let $a=(1,3,3,4)$ and $I,J\triangleleft k_a[w,x,y,z]$ be given by
    \[I = (w^2, wx+z, x^3+w^2yz),\quad J = (w^2, w^4+y).\]
    Then $I$ is weighted-homogeneous, but $J$ isn't, since $\deg w^4=4\neq 3=\deg y$.
\end{example}

\begin{lemma}[Equivalent definition of weighted-homogeneous ideals]\label{lem:equiv-def-wh-ideals}
    An ideal $I\triangleleft\kazn$ is weighted-homogeneous if and only if every element $f\in I$ can be written as
        \[f=\sum_{i=0}^{\deg f} f_i\]
    for unique $f_i\in\kazn_i\cap I$.
\end{lemma}

\begin{proof}
    \emph{This proof was originally going to be a simple reference to some pre-existing proof, but apparently everybody else has had the same idea -- the author couldn't find a reference which stated this fact and didn't leave its proof as an exercise to the reader.
    If this paper is good for nothing else, at least it might provide a reference for people looking for an easy-to-find proof of this fact.}

    \bigskip

    Note that, for any $f\in I$, writing\footnote{%
        When we say $f\cap\kazn_i$ we mean write $f=\sum f_j$, where each $f_j$ is a sum of monomials of degree $j$, and then define $f\cap\kazn_i=\{f_j\}\cap\kazn_i=f_i$.
    } $f_i=f\cap\kazn_i$ we have
    \[f=\sum_{i=0}^{\deg f}f_i\]
    with the $f_i$ uniquely determined by $f$.
    So the above condition is equivalent to requiring that $f\cap\kazn_i\in I$ for all $i$ (since this intersection is empty, and thus trivially in $I$, for $i>d$).
    Yet another way of phrasing this is that $I$ must satisfy $I = \sum_i (I\cap\kazn_i)$.

    If $I$ satisfies this above condition then for each $i$ we find $g_j^{(i)}\in(I\cap\kazn_i)$ such that
    \[
        (g_1^{(i)},\ldots,g_{n_i}^{(i)})=I\cap\kazn_i.
    \]
    By definition, each $g_j^{(i)}$ must be weighted-homogeneous of degree $i$.
    But then
    \[
        I = \sum_i (I\cap\kazn_i) = \sum_i (g_1^{(i)},\ldots,g_{n_i}^{(i)}) = (\cup_i\{g_1^{(i)},\ldots,g_{n_i}^{(i)}\})
    \]
    is a way of writing $I$ as being generated by weighted-homogeneous elements.
    So $I$ is weighted-homogeneous.

    For the other direction, we assume that $I$ is weighted-homogeneous, so $I=(g_1,\ldots,g_n)$ for some weighted-homogeneous $g_i$.
    Note that $\{g_i\}\subseteq\sum_i(I\cap\kazn_i)$.
    Thus
    \[
        \sum_i(I\cap\kazn_i) \subseteq I = (g_1,\ldots,g_n) \subseteq \sum_i(I\cap\kazn_i),
    \]
    and so these ideals must in fact be equal.
\end{proof}

\begin{lemma}\label{lem:wh-prime}
    A weighted-homogeneous ideal $I\triangleleft\kazn$ is prime if and only if, whenever $fg\in I$ for $f,g\in\kazn$ with $f,g$ both homogeneous, either $f\in I$ or $g\in I$.
    That is, when considering primality of the ideal, it is enough to check the usual definition on only the homogeneous elements of the ideal.
\end{lemma}

\begin{proof}
    We can simply appeal to the proof of \cite[Chapter~6,~Section~1,~Lemma~1.2]{Agrawal:uf}.
\end{proof}

The idea of saying that some polynomial in a weighted-homogeneous ideal vanishes at some point is still well defined, since every polynomial in the ideal is a sum of multiples of weighted-homogeneous polynomials.

    \subsection{Weighted projective varieties} 
    \label{sub:weighted_projective_varieties}
    
    Motivated by our previous discussion (defining what it means for a weighted-homogeneous polynomial to vanish at a point in weighted projective space) we can now define the idea of a \emph{weighted projective variety}, much in the same way as one would for straight projective space.

    \begin{definition}[Weighted projective varieties and their ideals]\label{defn:v-and-i}
        Let $I\triangleleft \kazn$ be a weighted-homogeneous ideal.
        Define the \emph{weighted projective variety (associated to $I$)} by
        \[\van(I) = \{p\in\pazn \mid f(p)=0\text{ for all }f\in I\}.\]

        Conversely, let $V\subseteq\pazn$ be a subset of weighted projective space.
        Define the \emph{ideal associated to $V$} by
        \[\ide(V) = \{f\in \kazn \mid f(p)=0\text{ for all }p\in V\text{ and }f\text{ is }a\text{-weighted-homogeneous}\}.\]

        We say that a subset $V\subseteq\pee(a)$ is a \emph{weighted projective variety} if it is of the form $\van(I)$ for some weighted-homogeneous ideal $I\triangleleft \kazn$.
        
        If we have two varieties $V\subseteq W$ then we say that $V$ is a \emph{subvariety} of $W$.
        A weighted projective variety is said to be \emph{irreducible} if it has no non-trivial decomposition into subvarieties, i.e. $V=V_1\cup V_2$ with $V_1,V_2\neq\varnothing,V$.
    \end{definition}

    It's important to point out that, although we call $\ide(V)$ the \emph{ideal} associated to $V$, we have yet to prove that it actually \emph{is} an ideal.
    We do this in \cref{lem:fuf}.

    \begin{note}
        As a slight notational quirk we write $\van\ide$ to mean the composition $\van\circ\ide$, since this looks a lot less messy.
        So instead of writing, for example, $\van(\ide(\van(I)))$, we write $\van\ide\van(I)$.
    \end{note}

    \begin{example}
        Here are some simple examples of $\van$ and $\ide$ using just the definitions that we have so far.
        We will develop and uncover some more advanced machinery and techniques later on in this section and the next.
        \begin{itemize}
            \item $\van(x_i)=\{|x_0:\ldots:x_n|\in\pazn \mid x_i=0\}$, and so $U_i=\pazn\setminus\van(x_i)$;
            \item Let $X = \van(x_1^{a_2}-x_2^{a_1})\subseteq\pathree$.
            Then if $x=|x_0:x_1:x_2|\in X$ we must have that $x_1^{a_2}=x_2^{a_1}$.
            Splitting this into two cases ($x_1,x_2\neq0$ and $x_1,x_2=0$) gives us
            \[
                X= \{|x_0:1:1|\}\cup\{|1:0:0|\}
            \]
            where we use that fact that $|x_0:x_1:x_2|=\left(1/x_1\right)^{1/a_1}\cdot|x_0:x_1:x_2|$ for the first case;
            \item What is $\ide\van(x_i^2)$?
            Well $x_i^2=0$ if and only if $x_i=0$, since $k$ is a field, so
            \[
                \van(x_i^2) = \{|x_0:\ldots:x_n|\in\pazn \mid x_i=0\} = \van(x_i).
            \]
            But the $x_j$ (for $j\neq i$) can take any value in $k$ (as long as they aren't all simultaneously zero), so the only polynomials $f\in\kazn$ that satisfy $f(x)=0$ for all $x\in\van(x_i^2)=\van(x_i)$ are those of the form $x_ig(x)$ for some polynomial $g\in\kazn$.
            That is,
            \[
                \ide\van(x_i^2) = \ide\van(x_i) = (x_i) \triangleleft \kazn.\qedhere
            \]
        \end{itemize}
    \end{example}

    We will later see that weighted projective varieties are in fact also projective varieties in the usual sense, and so they really are deserving of the name `varieties'.
    But then it doesn't seem too unreasonable to hope that we could define some sort of Zariski topology on weighted projective space with our weighted projective varieties, so we state and prove a lemma that lets us do so.

    \begin{lemma}\label{lem:topology-well-defined}
        Let $I,J\triangleleft\kazn$ be weighted-homogeneous ideals.
        Then
        \begin{enumerate}[(i)]
            \item $\van(I)\cup\van(J)=\van(IJ)$
            \item $\van(I)\cap\van(J)=\van(I+J)$
            \item $\varnothing=\van(\kazn)$;
            \item $\pazn=\van(\zi)$\qedhere
        \end{enumerate}
    \end{lemma}
    
    \begin{proof}
        \begin{enumerate}[(i)]
            \item An element of $IJ$ is by definition of the form $\sum^k f_ig_i$ for $f_i\in I$ and $g_i\in J$, which tells us that $\van(I)\cup\van(J)\subseteq\van(IJ)$, since if $x$ vanishes on either all of $I$ or all of $J$ (or maybe some of both) then it definitely vanishes on all elements of $IJ$.

            Conversly, $fg\in IJ$ for all $f\in I$ and $g\in J$.
            Since $\kazn$ is Noetherian, $I=(f_1,\ldots,f_k)$ and $J=(g_1,\ldots,g_l)$.
            So if $x$ vanishes on all of $IJ$ then in particular it vanishes on
            \[\{f_ig_j\mid 1\leqslant i\leqslant k, 1\leqslant j\leqslant l\} = \bigcup_{i=1}^k\{f_ig_j\mid 1\leqslant j\leqslant l\}.\]

            If $f_i(x)\neq0$ for some $i$ then we must have that $g_j(x)=0$ for all $x$ and so $x\in J$, and if not then $x\in I$.
            Thus $\van(IJ)\subseteq\van(I)\cup\van(J)$.
            \item If $x\in\van(I)\cap\van(J)$ then $f(x)=g(x)=0$ for all $f\in I$ and $g\in J$.
            Thus $(f+g)(x)=0$ for all $(f+g)\in I+J$ (which is all of $I+J$ by definition).

            Conversely, since $0\in I,J$, for all $f\in I$ and $g\in J$ we know that $f=f+0,g=0+g\in I+J$.
            So if $x\in\van(I+J)$ then it vanishes in particular on all of $I$ and all of $J$.
            \item The only `point' which vanishes on $x_i$ for all $i$ is $0$, but $0\not\in\pazn$, so
            \[\van(\kazn)\subseteq\van(x_0,\ldots,x_n)=\varnothing.\]
            \item Every point in $\pazn$ vanishes on the zero polynomial.\qedhere
        \end{enumerate}
    \end{proof}

    \begin{lemma}
        An arbitrary sum
        \[
            I = \sum_{\alpha\in\mathcal{A}}I_\alpha = \left\{\sum_{\beta\in\mathcal{B}}f_\beta \bigg\rvert f_\beta\in I_\beta\text{ and } \mathcal{B}\subset\mathcal{A}\text{ is finite}\right\}
        \]
        of weighted-homogeneous ideals is a weighted-homogeneous ideal.
    \end{lemma}
    
    \begin{proof}
        We first claim that an abitrary sum of ideals is an ideal.
        It is clear that $0\in I$, as well as $rf\in I$ for any $f\in I$ and $r\in\kazn$, so it remains only to show that $I$ is closed under finite sums.
        Let $f = \sum_{\beta\in\mathcal{B}}f_\beta$ and $g = \sum_{\beta\in\mathcal{C}}g_\beta$ be elements of $I$.
        We note that $f=\sum_{\beta\in\mathcal{D}}f'_\beta$, where $\mathcal{D}\supseteq\mathcal{B}$ is finite and
        \[
            f'_\beta =
            \begin{cases}
                f_\beta & \text{if }\beta\in\mathcal{B} ; \\
                0 & \text{if }\beta\in\mathcal{D}\setminus\mathcal{B} .
            \end{cases}
        \]
        So let $\mathcal{D}=\mathcal{B}\cup\mathcal{C}$.
        Then $f+g = \sum_{\beta\in\mathcal{D}}(f'_\beta+g'_\beta)\in I$, since $(f'_\beta+g'_\beta)\in I_\beta$ for each $\beta\in\mathcal{D}$.
        
        To show further that this sum is a weighted-homogeneous ideal we use \cref{lem:equiv-def-wh-ideals}.
        By definition, every element $f$ in the arbitrary sum $J=\sum_{I\in\mathcal{I}}I$ is a finite sum of elements in the summands.
        That is, there exist $I_1,\ldots,I_k\in\mathcal{I}$ and $f_i\in I_i$ such that
        \[f = \sum_{i=1}^k f_i.\]
        But each $I_i$ is weighted-homogeneous, and so each $f_i$ can be written as a sum of weighted-homogeneous elements $g_j^{(i)}$ where $\deg g_j^{(i)}=j$.
        Thus
        \[f = \sum_{j=1}^l\left(\sum_{i=1}^k g_j^{(i)}\right)\]
        is an expression for $f$ as a sum of weighted-homogeneous elements, and so $J$ is a weighted-homogeneous ideal.
    \end{proof}

    \begin{corollary}\label{cor:topology-well-defined-2}
        An arbitrary intersection of weighted projective varieties is a weighted projective variety:
        \[\bigcap_{I\in\mathcal{I}}\van(I)=\van\left(\sum_{I\in\mathcal{I}}I\right)=\van(J)\]
        where $\sum_{I\in\mathcal{I}}I=J\triangleleft\kazn$ is a weighted-homogeneous ideal.
    \end{corollary}

    With \cref{lem:topology-well-defined,cor:topology-well-defined-2} in hand, we can now make the definition that we would like to.

    \begin{definition}[Zariski topology]\label{defn:zariski-topology}
        The \emph{Zariski topology on $\pazn$} is given by defining the closed sets of $\pazn$ to be those of the form $\van(I)$ for some weighted-homogeneous ideal $I\triangleleft \kazn$, that is, the weighted projective varieties.
    \end{definition}

    One final thing to note before moving onto the Nullstellensatz is how we can use the construction of weighted projective space to understand these weighted projective varieties.
    The way that we define $f(p)=0$ for some $a$-weighted-homogeneous $f$ and point $p\in\pee(a)$ is really by requiring that $f(\hat{p})=0$, where $\hat{p}\in\aff^{n+1}\noz$ is a representative of $p$.
    We use the requirement of $f$ being $a$-weighted-homogeneous to ensure that this definition is well-defined under a change of representatives.

    So we can think of $\van(I)$ as a quotient of the affine `cone'\footnote{%
        It is only a cone though for $a_0,\ldots,a_n=1$, whence we recover straight projective space.
        The use of the word here is definitely meant to be understood in the sense of a loose meaning, but we will from now one speak of the affine cone anyway, to avoid having to constantly write cone in quotation marks.
    }:
    \begin{equation}\label{eq:affine-cone-thing}
        \van(I) = \frac{\vanaff(I)\noz}{\gm} \subseteq \frac{\aff^{n+1}\noz}{\gm}
    \end{equation}
    where $\vanaff(I)=\{x\in\aff^{n+1} \mid f(x)=0\text{ for all }f\in I\}$ and we consider $I\triangleleft k[x_0,\ldots,x_n]$ as an ideal in the usual polynomial ring (i.e. with all weights equal to $1$, though really this doesn't matter, since the weights affect only the graded structure of the ring).

    \begin{definition}[Affine cone]\label{defn:affine-cone-thing}
        Given $X=\van(I)$ for some weighted-homogeneous ideal $I\triangleleft\kazn$ we write $\hat{X}$ to mean $\vanaff(I)$, so that \cref{eq:affine-cone-thing} can be written as
        \[
            X = \frac{\hat{X}\cap(\aff^{n+1}\noz)}{\gm}.\qedhere
        \]
    \end{definition}

    Note that we don't simply write $\hat{X}\noz$ in \cref{defn:affine-cone-thing} since we don't know a priori that $0\in\hat{X}$, but we clear this up in \cref{lem:affine-ide-on-cone-same-as-ide-on-var}.


    \subsection{The weighted projective Nullstellensatz} 
    \label{sub:the_weighted_projective_nullstellensatz}

    \begin{lemma}[Five useful facts]\label{lem:fuf}
        Let $I,J\triangleleft \kazn$ be weighted-homogeneous ideals and let $V,W\subseteq\pazn$.
        Then with $\van$ and $\ide$ be defined as in \cref{defn:v-and-i} we have that
        \begin{enumerate}[(i)]
            \item $\ide(V)\subseteq \kazn$ is a radical weighted-homogeneous ideal;
            \item If $I\subseteq J$ then $\van(I)\supseteq\van(J)$;
            \item If $V\subseteq W$ then $\ide(V)\supseteq\ide(W)$;
            \item $I\subseteq\ide\van(I)$;
            \item\label{item:fuf-van-ide-X-is-X} $\van(I)=\van\ide\van(I)$.\qedhere
        \end{enumerate}
    \end{lemma}

    \begin{proof}
        \begin{enumerate}[(i)]
            \item The fact that $\ide(V)$ is an ideal is reasonably straightforward: if $f$ and $g$ vanish at all points of $V$ then so too does $f+g$, as does $fh$ for any other polynomial $h$.

            Since $\kazn$ is Noetherian we know then that $\ide(V)=(f_1,\ldots,f_m)$ for some $f_i\in\ide(V)$.
            But by definition, $f\in\ide(V)$ means that $f$ must be $a$-weighted-homogeneous.
            Thus $\ide(V)$ is generated by weighted-homogeneous elements, and so is an $a$-weighted-homogeneous ideal.

            Finally, say that $f^k\in\ide(V)$.
            Then $0=(f^k)(p)=f(p)^k$ and so $f(p)=0$ since $\kazn$ is an integral domain.
            So $f\in\ide(V)$ and hence $\ide(V)$ is radical.
            \item\label{item:fuf-van-reversing} If $p\in\van(J)$ then $f(p)=0$ for all $f\in J$, and thus for all $f\in I$, so $p\in\van(I)$.
            \item\label{item:fuf-ide-reversing} If some polynomial vanishes at all points in $W$ then it vanishes in particular at all points in $V\subseteq W$.
            \item\label{item:fuf-ideal-contained-in-IV} Let $f\in I$, so that by definition $f(x)=0$ for all $x\in\van(I)$ and thus \mbox{$f\in\ide\van(I)$}.
            \item Let $x\in\van(I)$, then $f(x)=0$ for all $f\in\ide\van(I)$, and so $x\in\van\ide\van(I)$.
            Conversely, use (\ref{item:fuf-van-reversing}) with (\ref{item:fuf-ideal-contained-in-IV}) to get that $\van\ide\van(I)\subseteq\van(I)$.\qedhere
        \end{enumerate}
    \end{proof}

    \begin{definition}[Relevant ideals]\label{defn:relevant-ideals}
        An ideal $I\triangleleft\kazn$ is said to be \emph{relevant} if it satisfies the following two conditions:
        \begin{enumerate}[(i)]
            \item it is strictly contained inside the irrelevant ideal\footnote{%
                For more on the irrelevant ideal (and graded rings in general) see \cite[Chapter~6]{Agrawal:uf}.
            } $(x_0,\ldots,x_n)$;
            \item $\van(I)\neq\varnothing$.\qedhere
        \end{enumerate}
    \end{definition}

    \begin{note}
        If $I$ is weighted-homogeneous then the first condition in \cref{defn:relevant-ideals} is actually redundant: it is enough to simply ask that $\van(I)\neq\varnothing$, since this implies that $I$ is strictly contained inside the irrelevant ideal.\footnote{Thanks to Christopher R. Miller for pointing this out and correcting previous versions of this definition, as well as providing some interesting examples.}
        To see this, we can use the affine Nullstellensatz and Equation~(\ref{eq:affine-cone-thing}) to get that
        \begin{equation*}
            (x_0,\ldots,x_n)\subseteq I \implies \vanaff(I)\subseteq\{0\} \iff \van(I)=\varnothing.
        \end{equation*}
        So $\van(I)\neq\varnothing$ implies that $(x_0,\ldots,x_n)\not\subseteq I$.
        In particular, $(x_0,\ldots,x_n)\neq I$.
        But if $I$ is weighted-homogeneous then $I\subseteq(x_0,\ldots,x_n)$, since all of its generators must be of the form $\sum a_ix_i$.
        Thus if $\van(I)\neq\varnothing$ then $I$ is strictly contained inside the irrelevant ideal.
    \end{note}
    
    Just as there are multiple equivalent definitions for an ideal to be homogeneous, there are multiple equivalent definitions for an ideal to be relevant.
    Sometimes one is more useful than the other, so we list four equivalent conditions here, and when we use `the' definition of a relevant ideal in a proof or suchlike we mean any one of the following conditions.

    \begin{lemma}\label{lem:relevant-other-defs}
        Let $I\triangleleft\kazn$ be a weighted-homogeneous ideal.
        Then the following are equivalent:
        \begin{enumerate}[(i)]
            \item $I$ is relevant;
            \item $I$ is strictly contained inside $\kazn$ and is not equal to the irrelevant ideal;
            \item $(x_0,\ldots,x_n)\not\subseteq\rad(I)$.\qedhere
        \end{enumerate}
    \end{lemma}

    \begin{proof}
        \begin{description}
            \item[$(i)\iff(ii)$:] If $I$ is strictly contained inside the irrelevant ideal then it is clearly strictly contained inside the whole ring.

            Conversely, if $I$ is strictly contained inside the whole ring then it cannot have any constant generators, since $(\theta)=\kazn$ for any $\theta\in k$.
            Thus $I$ is strictly contained inside the irrelevant ideal (since by assumption it is not equal to it).
            \item[$(i)\iff(iii)$:] By the the affine Nullstellensatz we know that $\ideaff\vanaff(I)=\rad(I)$, and so using \footnote{%
                Or if we are happy with the fact that the irrelevant ideal is radical then we could just take the $\rad$ of both sides of $(x_0,\ldots,x_n)\subseteq I$.
            } \cref{lem:fuf} (\ref{item:fuf-ide-reversing}) gives
            \begin{align*}
                \van(I)=\varnothing &\iff \vanaff(I)\subseteq\zi \\&\iff \ideaff(\zi)\subseteq\ideaff\vanaff(I) \iff (x_0,\ldots,x_n)\subseteq\rad(I).\qedhere
            \end{align*}
        \end{description}
    \end{proof}

    \begin{definition}[Maximal weighted-homogeneous ideals]\label{defn:max-wh-ideals}
        An ideal $I\triangleleft\kazn$ is said to be a \emph{maximal weighted-homogeneous ideal} if it is relevant and maximal amongst relevant weighted-homogeneous ideals.
        That is, if $J\triangleleft\kazn$ is a weighted-homogeneous ideal such that $I\subsetneq J$, then $J$ is irrelevant (so $J=(x_0,\ldots,x_n)$ or $J=\kazn$).
    \end{definition}

    With all of these definitions out of the way, we now start our journey towards the weighted projective Nullstellensatz.
    We do this by proving a few technical lemmas, and then the Nullstellensatz drops out quite easily and naturally from them.
    Really then, the way to understand this train of thought is to read \cref{thm:weighted-proj-null} \emph{first} and then come back to \cref{lem:radical-wh,cor:ideaff-cone-is-wh,lem:affine-ide-on-cone-same-as-ide-on-var}, otherwise it might seem like the lemmas are pulled from thin air.

    \begin{lemma}\label{lem:radical-wh}
        Let $I\triangleleft\kazn$ be a weighted-homogeneous ideal.
        Then $\rad(I)\triangleleft\kazn$ is also a weighted-homogeneous ideal.
    \end{lemma}

    \begin{proof}
        Let $f\in\rad(I)$, so that $f^k\in I$ for some $k\in\nn$.
        Write $d=\deg f$ and let
        \[
            f_i=f\cap\kazn_i\quad\text{for }0\leqslant i\leqslant d
        \]
        (since for $i>d$ this intersection will be empty).
        Then by \cref{lem:equiv-def-wh-ideals} it is enough to show that $f_i\in\rad(I)$ for all $0\leqslant i\leqslant d$, since the $f_i$ are uniquely determined by $f$.

        We look first at $f_d$.
        Because $f_d^k=f^k\cap\kazn_d$ (since it is the only term of high enough degree) and $I$ is weighted-homogeneous, we must have that $f_d^k\in I$, and so $f_d\in\rad(I)$.
        But then $f-f_d\in\rad(I)$ is a polynomial of strictly smaller degree with homogeneous components $f_0,\ldots,f_{d-1}$, thus $(f-f_d)^{k'}\in I$ for some $k'\in\nn$, so we repeat the above process with $f_{d-1}$ to show that $f_{d-1}^{k'}\in I$, and thus $f_{d-1}\in\rad(I)$.

        After repeating this finitely many times (since the total degree strictly decreases each time) we have that $f_i\in\rad(I)$ for all $0\leqslant i\leqslant d$.
    \end{proof}

    \begin{corollary}\label{cor:ideaff-cone-is-wh}
        Let $I\triangleleft\kazn$ be a weighted-homogeneous ideal.
        Then $\ideaff\vanaff(I)$ is also a weighted-homogeneous ideal.
    \end{corollary}

    \begin{proof}
        The affine Nullstellensatz tells us that $\ideaff\vanaff(I)=\rad(I)$, so by \cref{lem:radical-wh} we are done.
    \end{proof}

    \begin{lemma}\label{lem:max-wh-are-radical}
        Let $I\triangleleft\kazn$ be a maximal weighted-homogeneous ideal.
        Then $I$ is radical.
    \end{lemma}

    \begin{proof}
        By \cref{lem:relevant-other-defs} we know that $(x_0,\ldots,x_n)\not\subseteq\rad(I)$.
        Thus $\van(\rad(I))\not\subseteq\varnothing$, and $\van(\rad(I))\neq\varnothing$.
        So $\rad(I)$ is relevant.
        Further, $\rad(I)$ is a weighted-homogeneous ideal, by \cref{lem:radical-wh}.
        We also know that $I\subseteq\rad(I)$, but if $I$ is a proper subset of $\rad(I)$ then this gives us our contradiction, since $I$ is maximal amongst radical weighted-homogeneous ideals.
        Hence $I=\rad(I)$.
    \end{proof}

    The next lemma (\cref{lem:affine-ide-on-cone-same-as-ide-on-var}) is really the key to the weighted projective Nullstellensatz, but it is largely just technicalities and abstract faff, so we state and prove it as a separate lemma just to make the statement and proof of \cref{thm:weighted-proj-null} a bit more slick.

    \begin{lemma}\label{lem:affine-ide-on-cone-same-as-ide-on-var}
        Let $I\triangleleft \kazn$ be a weighted-homogeneous relevant ideal and $X=\van(I)$.
        Then
        \[
            \ide(X) = \ideaff(\hat{X}).\qedhere
        \]
    \end{lemma}

    \begin{proof}
        Let $f$ be a generator of $\ideaff(\hat{X})$ (of which there are finitely many, since $k[x_0,\ldots,x_n]$ is Noetherian).
        By \cref{cor:ideaff-cone-is-wh} we know that $f$ is $a$-weighted-homogeneous.
        Also $f(\hat{x})=0$ for all $\hat{x}\in\hat{X}$, and so in particular $f(\hat{x})=0$ for all $\hat{x}\in\hat{X}\noz$.
        Combining these two facts we see that $f\in\ide(X)$.
        Then, since all the generators of $\ideaff(\hat{X})$ are in $\ide(X)$,
        \[
            \ideaff(\hat{X})\subseteq\ide(X).
        \]
        
        For the other inclusion we use the fact that $I$ is relevant, and so $X=\van(I)\neq\varnothing$.
        We also know that $X=\van\ide(X)$, thus $\van\ide(X)\neq\varnothing$.
        Hence $\ide(X)$ is also relevant.
        So there are no constant polynomials in $\ide(X)$, since otherwise $\ide(X)$ would be the whole of $\kazn$, contradicting the fact that it is relevant.
        Hence if $f\in\ide(X)$ then $f(0)=0$.

        Also, if $f\in\ide(X)$ then $f(x)=0$ for all $x\in X$, i.e. $f(\hat{x})=0$ for all representatives $\hat{x}\in\hat{X}\noz$.
        So $f\in\ideaff(\hat{X}\noz)$.
        But since $f(0)=0$ as well, $f\in\ideaff(\hat{X})$, hence
        \[
            \ide(X)\subseteq\ideaff(\hat{X}).\qedhere
        \]
    \end{proof}

    \begin{theorem}[Weighted projective Nullstellensatz]\label{thm:weighted-proj-null}
        Let $I\triangleleft \kazn$ be a weighted-homogeneous relevant ideal.
        Then
        \[
            \ide\van(I)=\rad(I).\qedhere
        \]
    \end{theorem}

    \begin{proof}
        Write $X=\van(I)$.
        Then \cref{lem:affine-ide-on-cone-same-as-ide-on-var} tells us that
        \[
            \ide\van(I) = \ide(X) = \ideaff(\hat{X}) = \rad(I).\qedhere
        \]
    \end{proof}

    \begin{corollary}[Applied weighted projective Nullstellensatz]\label{cor:applied-wp-null}
        The maps $\van$ and $\ide$ give us an inclusion reversing bijection between weighted projective varieties and radical weighted-homogeneous relevant ideals:
        \[
            \begin{array}{rcl}
                \underbrace{\left\{
                \begin{array}{c}
                    \text{radical w.h. relevant ideals}\\
                    I\triangleleft\kazn
                \end{array}
                \right\}}_{\ide\van\text{ acts as the identity here}}
                &\begin{array}{c}
                    \overset{\van}{\longrightarrow}\\
                    \overset{\ide}{\longleftarrow}
                \end{array}
                &\underbrace{\left\{
                \begin{array}{c}
                    \text{weighted projective varieties}\\
                    \varnothing\neq X=\van(I)\subseteq\pazn
                \end{array}
                \right\}}_{\van\ide\text{ acts as the identity here}}\\[3em]
                I\subseteq J &\implies &\van(J)\subseteq\van(I)\\
                \ide(Y)\subseteq \ide(X) &\Longleftarrow &X\subseteq Y.
            \end{array}
        \]

        Further, under this bijection, prime weighted-homogeneous ideals correspond to irreducible varieties, and maximal weighted-homogeneous ideals\footnote{Recall \cref{defn:max-wh-ideals}.} to points.
    \end{corollary}

    \begin{proof}
        The first part of the statement follows directly from \cref{thm:weighted-proj-null}.
        
        \bigskip
        
        Next we consider the statement about prime ideals and irreducible varieties.
        Assume first that $I\triangleleft\kazn$ is a prime weighted-homogeneous relevant ideal.
        Let
        \[
            X=\van(I)=\van(I_1)\cup\van(I_2)=X_1\cup X_2
        \]
        be a decomposition of $X$, and assume that $I_i$ is radical (since $\van(I_i)=\van(\rad(I_i))$ by taking $\van$ of both sides of \cref{thm:weighted-proj-null}).
        We want to show that $X_1=\varnothing$ or $X_1=X$.
        Since $X_i\subseteq X$ we know that $\van(X)\subseteq\van(X_i)$, i.e. $I\subseteq I_i$.
        We also know that $\van(I_1)\cup\van(I_2)=\van(I_1 I_2)$, and thus $I=I_1 I_2$ (since the product of radical ideals is again radical).
        But $I$ is prime, and so, since (trivially by the above equality) $I_1 I_2\subseteq I$, we must have $I_i\subseteq I$ (without loss of generality assume that $I_1\subseteq I$).
        So $I_1\subseteq I\subseteq I_1$, hence $I_1=I$ and $X_1=X$.

        For the other direction, assume that $X=\van(I)$ is an irreducible weighted projective variety and let $fg\in I$ be such that $f,g\not\in I$.
        Then $X\subseteq\van(f)\cup\van(g)$, and so
        \[
            X = (X\cap\van(f)) \cup (X\cap\van(g))
        \]
        is a non-trivial decomposition, because $f,g\not\in\ I=\ide(X)$, and thus $X\cap\van(f)\subsetneq X$ (and similarly for $g$).

        \bigskip

        Finally, the statement concerning maximal weighted-homogeneous ideals and points.
        Let $I\triangleleft\kazn$ be a maximal weighted-homogeneous ideal in the sense of \cref{defn:max-wh-ideals}.
        In particular then, $I$ is relevant, and thus $\van(I)\neq\varnothing$, so let $x\in\van(I)$.
        Then $\ide\van(I)\subseteq\ide(\{x\})$, and using \cref{lem:max-wh-are-radical} we see that $\ide\van(I)=\rad(I)=I$.
        Since $J$ is maximal weighted-homogeneous ideal, we know that either $\ide(\{x\})$ is irrelevant, or $\ide(\{x\})=I$.

        Say that $\ide(\{x\})$ is irrelevant, then $\{x\}=\van\ide(\{x\})=\varnothing$ is a contradiction.
        So we must have that $\ide(\{x\})=I$.
        Thus $\van(I)=\van\ide(\{x\})=\{x\}$ is a point.

        Now let $p=|p_0:\ldots:p_n|\in\pazn$ be a point and define the ideal $J \triangleleft \kazn$ by
        \[
            J = \bigcup_{\substack{k=0\\k\neq i}}^n \left(p_i^{a_k}x_k^{a_i}-p_k^{a_i}x_i^{a_k}\right)
        \]
        so that $J$ is weighted-homogeneous with $J\subseteq\ide(\{p\})$.
        Thus $\ide(\{p\})\neq\varnothing$, and hence $\ide(\{p\})$ is relevant.
        Let $\mathfrak{a}\triangleleft\kazn$ be a relevant weighted-homogeneous ideal with $\ide(\{p\})\subseteq\mathfrak{a}$.
        Then $\varnothing\subsetneq\van(\mathfrak{a})\subseteq\van\ide(\{p\})=\{p\}$, and so $\van(\mathfrak{a})=\{p\}$.
        This then tells us that $\ide(\{p\})=\rad(\mathfrak{a})$.
        But $\mathfrak{a}\subseteq\rad(\mathfrak{a})=\ide(\{p\})\subseteq\mathfrak{a}$.
        Hence $\mathfrak{a}=\ide(\{p\})$, and $\ide(\{p\})$ is a maximal weighted-homogeneous ideal.
    \end{proof}

    \begin{note}
        Since the bijection in \cref{cor:applied-wp-null} is only between \emph{non-empty} weighted projective varieties and radical weighted-homogeneous \emph{relevant} ideals, most of the theorems that we cover from now on concerning a weighted projective variety $X$ will include the hypothesis that $X$ is non-empty.
        Usually these theorems will be trivially true if $X=\varnothing$, but it is important to note that \textbf{the proofs we give assume (if stated) that $X=\van(I)$ is non-empty and thus $I$ is relevant.}
    \end{note}


    \subsection{Coordinate rings} 
    \label{sub:coordinate_rings}

    We define the idea of the coordinate ring of a weighted projective variety in this section, but it won't be until \cref{sub:explaining_} that we have a proper understand of it, or even much of a use.

    \begin{definition}[Weighted-homogeneous coordinate rings]
        Let $X=\van(I)$ be a non-empty weighted projective variety.
        Then define the \emph{weighted-homogeneous coordinate ring of $X$} to be
        \[
            S(X) = \frac{\kazn}{\ide(X)}.\qedhere
        \]
    \end{definition}

    If we write $A(Y)$ to mean the coordinate ring of an affine variety $Y$ then
    \[
        S(X) = \frac{\kazn}{\ide(X)} = \frac{\kzn}{\ideaff(\hat{X})} = A(\hat{X}).
    \]
    That is, the weighted-homogeneous coordinate ring of a non-empty weighted projective variety is simply the coordinate ring of its affine cone, but with the $a$-weighted grading.

    It's important to realise that we can't really think of the elements of $S(X)$ as polynomial functions on $X$, since `evaluating a function at a point in weighted projective space' is not a well-defined concept in general.
    They can be thought of as polynomial functions on the affine cone though, but this isn't always much use, since a morphism of affine cones doesn't necessarily descend everywhere to a morphism of weighted projective varieties (this will be covered more in \cref{sub:morphisms_between_varieties}).

    Equally important is the fact that isomorphic weighted projective varieties might have non-isomorphic weighted-homogeneous coordinate rings (again, covered in \cref{sub:morphisms_between_varieties}).
    


\section{An algebraic approach} 
\label{sec:algebra}

    \subsection{Explaining \texorpdfstring{$\proj$}{Proj} with the Nullstellensatz} 
    \label{sub:explaining_}

    \emph{For this section we take for granted a knowledge of graded rings (see \cite[Chapter~6]{Agrawal:uf}).}
    \emph{Many of the theorems in this section could be stated and proved for general graded rings $R$, but we tend to consider only the case when $R$ is also a finitely-generated $k$-algebra, to avoid straying too far into the world of schemes.}

    \emph{Recall \cref{ssub:list_of_notation_and_assumptions} -- when we say that a ring is graded we mean specifically $\zz^{\geqslant0}$-graded.}

    \bigskip

    The $\proj$ construction is defined for any graded ring $R$, but here we only really define it for specific types of graded rings\footnote{%
        Though for the sake of (a poor effort towards) completeness we will mention the more general definition afterwards.
    }, namely quotients of $\kazn$ by radical $a$-weighted-homogeneous relevant ideals $I\triangleleft\kazn$, that is, coordinate rings $S(\van(I))$.
    Note that requiring $I$ to be radical simply means that $k\subsetneq R$, so that $R$ is in a sense non-trivial.

    \begin{definition}[$\proj$ of a finitely-generated $k$-algebra]
        Let $R$ be a finitely-generated $k$-algebra of the form
        \begin{equation}\label{eqn:form-of-r}
            R = \kzny = \frac{\kazn}{I}
        \end{equation}
        where $y_i$ is the image of $x_i$ in the quotient by $I$ (so $\wt y_i=a_i=\wt x_i$), and $I$ is a radical $a$-weighted-homogeneous relevant ideal.
        We define \emph{the set}\footnote{%
            We make this distinction here because the natural definition of $\proj R$ is as a \emph{scheme}, i.e. with a topology and a structure sheaf.
            This doesn't matter too much here, since we try to stick to `hands-on' algebraic geometry as much as possible, but it is well worth bearing in mind.
            See \cite[Chapter~4.5]{Vakil:2015wa} for more.
        } $\proj R$ by
        \[
            \proj R = \{\prid\triangleleft R \mid \prid\text{ is an $a$-weighted-homogeneous prime ideal with }R_+\not\subseteq\prid\}
        \]
        where $R_+=(y_0,\ldots,y_n)$ denotes the irrelevant ideal.
    \end{definition}

    \begin{note}
        The more general definition of $\proj R$, for any graded ring $R$ with irrelevant ideal $R_+$, is where we define the set
        \[
            \proj R = \{\prid\triangleleft R \mid \prid\text{ is a homogeneous prime ideal with }R_+\not\subseteq\prid\}.\qedhere
        \]
    \end{note}

    \begin{example}[$\proj\kazn$]\label{ex:proj-kazn}
        Let's start off by looking at what seems like it might be the simplest example: when \mbox{$I=\zi$.}

        By \cref{cor:applied-wp-null} we know that weighted-homogeneous prime ideals $\prid\triangleleft\kazn$ correspond to irreducible varieties in $\pazn$, and that maximal weighted-homogeneous ideals correspond to points.
        So if we split up $\proj\kazn$ into maximal and prime-but-not-maximal weighted-homogeneous ideals and use this bijection coming from $\van$ and $\ide$ then we can consider $\pazn$ as a set of \emph{points in weighted projective space} and \emph{irreducible weighted projective varieties}:
        \begin{align*}
            \proj&\kazn =\\[0.5em]
            &\begin{array}{rcl}
                \underbrace{\{p\in\pazn\}}_{\van(\maid)\text{ where }\maid\text{ is maximal weighted-homogeneous}} & \sqcup & \underbrace{\{X\subseteq\pazn\}.}_{\van(\prid)\text{ where }\prid\text{ is prime-but-not-maximal weighted-homogeneous}}
            \end{array}
        \end{align*}

        So if we just consider maximal ideals then $\proj\kazn$ is simply $\pazn$, but when we throw in the other prime ideals as well we enrich the structure slightly: it is a set containing all the points of $\pazn$ \emph{as well as} all the irreducible weighted projective varieties in $\pazn$.
    \end{example}

    \begin{note}
        We could use this as an alternative definition of $\pazn$:
        \[
            \pazn=\proj\kazn.
        \]
        In fact, we might as well do so from now on, but when we speak of `points in $\pazn$' we still mean points in the old sense, and we still refer to the points of $\proj\kazn$ that correspond to varieties as `varieties'.

        This is an almost identical situation to when we define $\aff^n=\spec\kzn$.
        For more on this, see \cite[Sections~1.3,~1.6]{Reid:1995tu} or \cite[Section~3.2]{Vakil:2015wa}
    \end{note}

    So it seems like just considering maximal ideals will give us pretty much the whole picture of what $\proj R$ looks like, and the prime ideals will tell us what the varieties \emph{inside} $\proj R$ look like.
    We now make this rigorous.

    Let $R$ be a finitely-generated algebra as in \cref{eqn:form-of-r} and $\prid\triangleleft R$ an $a$-weighted-homogeneous prime ideal not containing $R_+=(y_0,\ldots,y_n)$.
    Then $\prid$ corresponds uniquely to the weighted-homogeneous prime ideal $\bar{\prid}\triangleleft\kazn$ with $I\subseteq\bar{\prid}$ such that $\prid=\bar{\prid}/I$.\footnote{%
        See \cite[Chapter~III,~Section~8,~Theorem~11]{Zariski:1958ue}.
    }
    Further, since $(y_0,\ldots,y_n)\not\subseteq\prid$ we know that $(x_0,\ldots,x_n)\not\subseteq\bar{\prid}$\footnote{%
        Proof by contradiction, using the fact that if we have ideals $\mathfrak{i}\subseteq\mathfrak{a}\subseteq\mathfrak{b}$ then $\mathfrak{a}/\mathfrak{i}\subseteq\mathfrak{b}/\mathfrak{i}$
    }, i.e. $\bar{\prid}$ is relevant (since $\bar{\prid}$ is prime and thus radical).
    The situation unfolds in the same way when $\maid\triangleleft R$ is a maximal weighted-homogeneous ideal, giving us a unique maximal weighted-homogeneous relevant ideal $\bar{\maid}\triangleleft\kazn$ such that $I\subseteq\bar{\maid}$ and $\maid=\bar{\maid}/I$.

    By \cref{cor:applied-wp-null}, $\bar{\maid}$ corresponds to a point $p_{\bar{\maid}}=\van(\bar{\maid})\in\pazn$, but $I\subseteq\bar{\maid}$ means that $p_{\bar{\maid}}\subseteq\van(I)$.
    Conversely, given any point in $\van(I)$ we see that it corresponds to a maximal weighted-homogeneous relevant ideal of $\kazn$ containing $I$, and thus to a maximal weighted-homogeneous ideal of $R$ \emph{not} containing $R_+$.
    Similarly, the inclusion-reversing bijection tells us that $\bar{\prid}$ corresponds to an irreducible weighted projective variety contained inside $\van(I)$, i.e. a subvariety of $\van(I)$, and vice versa.

    So we have the bijective correspondence

    \begin{equation}\label{eqn:proj-corres-1}
        \begin{array}{c}
            \{\prid\in\proj R\}\\
            \updownarrow\\
            \{X_{\bar{\prid}}\subseteq\van(I) \mid X_{\bar{\prid}}\text{ is an irreducible subvariety}\}
        \end{array}
    \end{equation}

    and, in particular,

    \begin{equation}\label{eqn:proj-corres-2}
        \begin{array}{c}
            \{\maid\in\proj R \mid \maid\triangleleft R\text{ is maximal weighted-homogeneous}\}\\
            \updownarrow\\
            \{p_{\bar{\maid}}\in\van(I) \mid p_{\bar{\maid}}\text{ is a point}\}.
        \end{array}
    \end{equation}

    \begin{theorem}\label{thm:what-proj-means-for-quotients}
        Let $I\triangleleft\kazn$ be a radical weighted-homogeneous relevant ideal.
        Then
        \[
            \proj\frac{\kazn}{I} = \van(I) \subseteq \pazn,
        \]
        but where $\proj R$ has this enriched structure of also containing points corresponding to all the irreducible subvarieties of $\van(I)$.

        Equivalently, let $X=\van(I)\subseteq\kazn$.
        Then
        \[
            \proj S(X) = \proj\frac{\kazn}{\ide(X)} = X \subseteq \pazn
        \]
        where we have the same enriched structure as above.
    \end{theorem}

    \begin{proof}
        This is just using \cref{eqn:proj-corres-1,eqn:proj-corres-2} in the same way as in \cref{ex:proj-kazn}.
    \end{proof}

    So in some sense\footnote{%
        The way of making this precise is to use category theory, which turns out to have very exciting applications to algebraic geometry as a whole.
        There is no quick introduction to this (at least, not that the author can find), but any good book on scheme theory should cover it, but \cite{Vakil:2015wa} is a particularly good text treating algebraic geometry \emph{after} covering a reasonable chunk of category theory.
        Failing that, it seems highly unlikely that \cite{Hartshorne:1977we} wouldn't cover it (and the author apologises for this infuriatingly vague reference).
    } $\proj$ and $S(\blank)$ are mutual inverses, i.e. $\proj S(X)) = X$ and $S(\proj R)) = R$.


    \subsection{Morphisms between varieties} 
    \label{sub:morphisms_between_varieties}

    \emph{The idea of morphisms is more complicated than it might sound at first.
        What we cover here is but a brief part of the whole story, as we take only what we need.
        Due to time, we (apologetically) might skim over some details, and this section is intended to be more of a motivation to read other sources than a complete guide.
        For more (and better) information, see \cite[Chapter~II,~Section~2]{Hartshorne:1977we}, where the real treatment is given using the language of schemes.}

    \emph{In summary: treat this section as a brief vacation from the usual rigour of mathematics to the land of allegory.}

    \bigskip

    We now have some notion of varieties, but as of yet have no rigorous idea of how we should define morphisms between them.
    Taking a cue from the affine and straight projective cases we think that a preliminary definition could be a map given by a polynomial in each coordinate.
    But before we make this definition, a thought occurs to us: we've just formalised the underlying algebraic structure of these geometric objects, and in the affine case a morphism of the coordinate rings induces a morphism of the varieties.
    So let's take this as a definition for the moment and see where we can get with it.

    \begin{definition}[Morphism of weighted projective varieties, attempt 1]\label{defn:morph-wpv-from-ring}
        Let $X\subset\pazn$ and $Y\subset\pee(b_0,\ldots,b_m)$ be weighted projective varieties, and $f\colon S(Y)\to S(X)$ a graded ring homomorphism.
        Write $S(Y)=k_{(b_0,\ldots,b_m)}[y_0,\ldots,y_m]/J$ and $S(X)=\kazn/I$, where $J=\ide(Y)$ and $I=\ide(X)$, so that $f(y_i)=f_i\in\kazn$.

        Then we define the \emph{morphism $F$ of weighted projective varieties} as
        \begin{align*}
            F_\#\colon X&\to Y \\
            x&\mapsto |f_0(x):\ldots:f_n(x)|.\qedhere
        \end{align*}
    \end{definition}

    We might (and should) worry about whether or not \cref{defn:morph-wpv-from-ring} is well defined.
    Is it always true that the $f_i$ never all vanish simultaneously?
    Is this map invariant under a different choice of representatives for the point $x\in\pazn$?
    Both of these questions can be answered satisfactorily, but we don't do so here, for reasons explained at the end of this section.
    Sweeping any and all problems of this sort under the proverbial rug, we march onwards.

    \begin{definition}[Isomorphism of weighted projective varieties]\label{defn:iso-wpv}
        Let $F\colon X\to Y$ be a morphism of weighted projective varieties.
        Then $F$ is an \emph{isomorphism} if there exists another morphism of weighted projective varieties $G\colon Y\to X$ such that $G\circ F=\id_X$ and $F\circ G=\id_Y$.
    \end{definition}

    But here we hit what seems like a problem: different embeddings of the same variety should definitely be isomorphic by any sensible definition, but we will see that different embeddings might not necessarily have isomorphic coordinate rings.
    For example, we will see in \cref{sub:truncation_of_graded_rings} that
    \[
         \pee^1 = \proj k[x,y] \cong \proj k[x,y]^{(2)} = \van(v^2-uw)\subset\pee^3
    \]
    but the former has coordinate ring $k[x,y]$, and the latter has coordinate ring $k[u,v,w]/(v^2-uw)$.
    Since the number of generators is different in each, the two definitely can't be isomorphic as graded rings.

    The way to solve this problem is to point out the following: \emph{not every morphism of varieties comes form a morphism of coordinate rings.}
    So \cref{defn:morph-wpv-from-ring} sounds great as a partial definition, i.e. that all maps of this form are indeed what we should call a morphism of varieties, but there are other maps that we should also call varieties.
    It turns out that listening to our original idea of polynomial maps would be sensible.

    \begin{definition}[Morphism of weighted projective varieties, attempt 2]
        Let $F\colon X\to Y$ be a map of weighted projective varieties, where $X\subset\pazn$ and $Y\subset\pee(b_0,\ldots,b_m)$.
        Then $F$ is a \emph{morphism} if it is a weighted-homogeneous polynomial in each coordinate.
        That is,
        \begin{align*}
            F\colon X&\to Y \\
            |x_0:\ldots:x_n| &\mapsto |F_0(x_0,\ldots,x_n):\ldots:F_m(x_0,\ldots,x_n)|
        \end{align*}
        where $F_i\in\kazn$ is weighted-homogeneous.
    \end{definition}

    Now we think about \cref{defn:iso-wpv}.
    It still makes sense as a definition, but there is another way that we could maybe define an isomorphism, using the algebraic structure again.
    If we had some bijection between prime (and maximal) ideals of the two coordinate rings of our varieties that preserved enough information, such as inclusion, then the varieties should be isomorphic, since the $\proj$ of their coordinate rings will be `the same', in a sense.
    This happens as an example in \cref{thm:proj-truncation-iso}.

    \bigskip

    It may seem like we are skirting around the issue of settling on a specific definition, and that's because we are.
    Really, the only times we will talk about morphisms in this paper is when we are talking about isomorphisms, and then we will come across just two cases:
    \begin{enumerate}[(i)]
        \item the underlying graded rings are isomorphic;
        \item one of the underlying graded rings is a truncation of the other.
    \end{enumerate}
    In case (i) it is very fair that, if two rings are isomorphic, then their $\proj$ should be isomorphic.
    In case (ii), \cref{thm:proj-truncation-iso} will apply.

    Unfortunately, since not much more theory than this is needed in this paper, not much more theory than this is covered in this paper.
    The real story is one of \emph{morphisms of schemes}, and once again the author recommends the invaluable resource that is \cite[Chapter~II]{Hartshorne:1977we}.
    

    \subsection{Truncation of graded rings} 
    \label{sub:truncation_of_graded_rings}

    We know that a variety, in the sense we've been describing them so far, depends on its ambient space, and from our experiences with affine varieties we expect that we might be able to embed the same variety in different weighted projective spaces.
    A classical example is that of a Veronese embedding of a variety from $\pee^n\to\pee^m$ for some specific $m\geqslant n$.
    In fact, this is the example that we study in \cref{ex:veronese-grading-thing}.

    So since \cref{thm:what-proj-means-for-quotients} gave us a way of converting between algebra and geometry, it seems like we should be able to find some process that we can apply to our coordinate rings that corresponds to an embedding of the associated varieties.

    \begin{definition}[Truncation of a graded ring]\label{defn:truncation-of-graded-ring}
        Let $R=\oplus_{i\geqslant0}R_i$ be a graded ring.
        Define $R^{(d)}$, the \emph{$d$-th truncation of $R$}, by
        \[
            R^{(d)} = \bigoplus_{i\geqslant0}R_{di}.
        \]
        So $R^{(d)}$ is also a graded ring, with grading given by $i$.
        That is, an element that has degree $di$ in $R$ has degree $i$ in $R^{(d)}$.\footnote{%
            The choice of which grading to use in the truncated ring (either the one that we use here, or simply keeping the grading the same) varies from author to author.
            As long as you pick one and stick with it doesn't (at least, not as far as this author knows) matter.
        }
    \end{definition}

    \begin{example}\label{ex:veronese-grading-thing}
        Let $R=k[x,y]$ with the usual grading (i.e. $\wt x,y=1$).
        Then
        \[
            R^{(2)} = \bigoplus_{i\geqslant0}R_{2i} = \bigoplus_{i\geqslant0}\{f\in k[x,y] \mid \deg f=2i\}.
        \]
        We note then that all polynomials in $R$ of even degree (which are exactly those that are in $R^{(d)}$) are generated by $x^2,xy,y^2$.
        Thus
        \[
            k[x,y]^{(2)} = k[x^2,xy,y^2].
        \]

        Now we know that $\proj k[x,y]=\pee(1,1)=\pee^1$, but what is $\proj k[x,y]^{(2)}$?
        First, let's write the latter in a form that we know how to deal with:
        \[
            k[x,y]^{(2)} = k[x^2,xy,y^2] \cong \frac{k[u,v,w]}{(uw-v^2)}.
        \]
        By \cref{defn:truncation-of-graded-ring} we have that $\deg x^2,xy,y^2=1$, and so taking $\wt u,v,w=1$ gives us an isomorphism of graded rings.
        In particular, it's important to note that $R\not\cong R^{(2)}$ here, and in general $R\not\cong R^{(d)}$.

        But now we can use \cref{thm:what-proj-means-for-quotients}:
        \[
            \proj\frac{k[u,v,w]}{(uw-v^2)} = \van(uw-v^2)\subseteq\pee(1,1,1) = \pee^2.
        \]
        This is exactly the degree-$2$ Veronese embedding of $\pee^1\hookrightarrow\pee^2$.
        So
        \[
            \pee^1 = \proj k[x,y] \cong \proj k[x,y]^{(2)} = \nu_2(\pee^1)\subseteq\pee^2
        \]
        and hence $\proj R\cong\proj R^{(2)}$ when $R=k[x,y]$, and in fact this truncation corresponds exactly to the degree-$2$ Veronese embedding.
    \end{example}

    It turns out that the above example is more than just a lucky coincidence.
    We have two claims:
    \begin{enumerate}[(i)]
        \item $\proj R\cong\proj R^{(d)}$ for any graded ring $R$ and $d\in\nn$;
        \item $\proj\kzn^{(d)}$ (so with $\wt x_i=1$) corresponds to the degree-$d$ Veronese embedding $\pee^n\hookrightarrow\pee^{\binom{n+d}{d} - 1}$.
    \end{enumerate}

    The second claim is slightly off-topic in a sense, since it is a fact concerning only straight projective space, but we can formulate it to deal with straight projective varieties too.
    It is a very nice example, but we unfortunately don't have the time to delve into it here any further.

    The first claim is one that we can apply to our studies of weighted projective varieties, and so we study it now.

        \subsubsection{The first claim} 
        \label{ssub:the_first_claim}

        \begin{theorem}\label{thm:proj-truncation-iso}
            Let $R$ be a graded ring and $d\in\nn$.
            Then
            \[\proj R\cong \proj R^{(d)}.\qedhere\]
        \end{theorem}

        \begin{proof}
            \emph{All that we really use in this paper is that $\proj R=\proj R^{(d)}$ \emph{as sets}.
            Since really they are endowed with so much more structure than we are covering here (namely, their structure as \emph{schemes}) the actual proof of this statement uses some ideas and definitions that we haven't mentioned, and don't intend to, for the sake of time.
            One important thing that we don't show is that $R^{(d)}[f^{-d}]=R[f^{-1}]^{(d)}$, and thus that in particular the degree-$0$ graded parts are equal.
            Rather than picking apart a standard proof and discarding the bits that we won't use, and hence potentially taking things out of context and missing the bigger picture, we provide here a sketch proof with commentary, lifted largely from \cite[Exercise~9.5]{Eisenbud:1995tm}.
            All of the statements and theorems concerning integral extensions can be found in \cite[Chapter~4.4]{Eisenbud:1995tm}.}

            \emph{For a full proof of this statement see \cite[Proposition~5.5.2]{Tevelev:jdT35_ao}, for example.
            And, of course, this (and so much more) is covered in \cite[Proposition~(2.4.7)]{Grothendieck:1961tr}.}

            \bigskip

            First of all we note that there is an injective graded ring homomorphism $R^{(d)}\hookrightarrow R$ corresponding to inclusion, but this is \emph{not} usually an isomorphism.
            So rather than looking at the underlying rings themselves, we look instead at the structure of their prime ideals.

            The above ring extension $R^{(d)}\hookrightarrow R$ is in fact integral, and thus every prime ideal of $R^{(d)}$ is given by the restriction of some prime ideal in $R$ to $R^{(d)}$, but this is in general \emph{not} be a bijective correspondence.
            However, it can be shown\footnote{%
                This is \emph{not} the author saying that this proof is left to the reader; this \emph{is} the author saying that the details are best left to a text on commutative algebra, and not an expository text on weighted projective spaces written by yours truly.
            } that when we consider only weighted-homogeneous prime ideals, we do in fact end up with a bijection $\prid\mapsto\prid\cap R^{(d)}$ from weighted-homogeneous prime ideals of $R$ to those in $R^{(d)}$.
        \end{proof}

        We can use \cref{thm:proj-truncation-iso} to simplify certain weighted projective spaces or varieties in two different ways:
        \begin{enumerate}[(a)]
            \item reduce $\pazn$ to a \emph{well-formed} weighted projective space $\pee(a_0',\ldots,a_n')$;
            \item embed $\pazn\hookrightarrow\pee^N$ for some large enough $N$ (`straighten out' $\pazn$).
        \end{enumerate}
        Both of these terms will be defined and explained next, before we approach an explicit example and work through it as best we can in \cref{sub:a_worked_example}.
        We start with well-formed weighted projective spaces in \cref{ssub:well_formed_weighted_projective_spaces} and then deal with `straightening out' $\pazn$ in \cref{subsub:embedding_}.
        Really, the second is sort of a specific case of the first, since straight projective space is a well-formed weighted projective space, but in another sense it is different entirely, since we want to end up in a specific weighted projective space: straight projective space.

        \subsubsection{Well-formed weighted projective spaces} 
        \label{ssub:well_formed_weighted_projective_spaces}

        \begin{definition}[Well-formed weights]
            We say that a weight $a=(a_0,\ldots,a_n)$ is \emph{well-formed} if any $n-1$ of the $a_i$ are coprime.
            That is,
            \[
                \gcd(a_0,\ldots,\remove{a_i},\ldots,a_n) = 1
            \]
            for all $0\leqslant i\leqslant n$.
        \end{definition}

        \begin{definition}[Well-formed weighted projective space]
            The weighted projective space $\pazn$ is said to be \emph{well-formed} if the weight $(a_0,\ldots,a_n)$ is well-formed.
        \end{definition}

        \begin{theorem}\label{thm:can-always-well-form}
            Given some weight $a$ there is a well-formed weight $a'$ such that $\pee(a)\cong\pee(a')$.
            That is, any weighted projective space $\pazn$ is isomorphic to a well-formed weighted projective space $\padzn$.
        \end{theorem}

        \begin{proof}
            Let $R=\kazn$, so that $\proj R=\pazn$.
            Since we might as well assume that $a$ is not already well-formed (otherwise the proof is trivial) we have two possible cases:
            \begin{enumerate}[(i)]
                \item there exists some common factor $d$ of all the $a_i$;
                \item $a_0,\ldots,a_n$ have no common factor, but there is some $j$ such that $a_0,\ldots,\remove{a_j},\ldots,a_n$ have common factor $d$ which is coprime to $a_j$.
            \end{enumerate}
            We use \cref{thm:proj-truncation-iso} for both cases.

            \bigskip

            In case (i) we have, for all $0\leqslant i\leqslant n$, that $d\mid a_i$, and thus $x_i\in R_{dk}$ for all $k$.
            Hence $R^{(d)}=k_{a/d}[x_0,\ldots,x_n]$.
            Thus
            \[
                \pazn =\proj\kazn \cong \proj k_{a/d}[x_0,\ldots,x_n] = \pee\left(\frac{a_0}{d},\ldots,\frac{a_n}{d}\right).
            \]

            \bigskip

            In case (ii) we see that, since $d$ is coprime to $a_j$, the only $a_j$ term that will appear in $R_{dk}$ is $a_j^d$ and its powers.
            So $R^{(d)}=k_{a/d}[x_0,\ldots,x_j^d,\ldots,x_n]$, and thus
            \[
                \pazn =\proj\kazn \cong \proj k_{a/d}[x_0,\ldots,x_j^d,\ldots,x_n] = \pee\left(\frac{a_0}{d},\ldots,a_j,\ldots,\frac{a_n}{d}\right).
            \]

            \bigskip

            So, given some weight $a=(a_0,\ldots,a_n)$, if they all share some common factor $d$ then we can use case (i) to divide all the $a_i$ through by $d$.
            We can repeat this until the $\gcd(a_0,\ldots,a_n)=1$.
            Then, if any $n-1$ of the $a_i$ have some common factor $d'$, we can use case (ii) to divide $a_0,\ldots,\remove{a_j},\ldots,a_n$ through by $d'$ until $\gcd(a_0,\ldots,\remove{a_j},\ldots,a_n)=1$.
        \end{proof}

        So in light of \cref{thm:can-always-well-form} there is usually no loss in generality in assuming that a weighted projective space is well-formed, unless we care about the specific embedding.
        We will \emph{not} always assume that $\pazn$ is always well-formed, and if we ever do then we will explicitly say so.

        \begin{example}
            Here are two particularly nice cases, the second of which will crop up again in \cref{sub:a_worked_example}:
            \begin{itemize}
                \item $\pee(a,b)\cong\pee(1,b)\cong\pee(1,1)=\pee^1$ for any $a,b\in\nn$;
                \item $\pee(ab,bc,ca)\cong\pee(b,bc,c)\cong\pee(1,c,c)\cong\pee(1,1,1)=\pee^2$ for any $a,b,c\in\nn$ (assumed to be coprime, without loss of generality).\qedhere
            \end{itemize}
        \end{example}
        

        \subsubsection{Embedding \texorpdfstring{$\pee(a_0,\ldots,a_n)$}{P(a0,...,an)} into \texorpdfstring{$\pee^N$}{PN}} 
        \label{subsub:embedding_}

        \emph{A good source that also covers the same material as in this section, and where the author first read most of this, is \cite[Section~5.5]{Tevelev:jdT35_ao}.}

        \bigskip

        Now we take a look at using \cref{thm:proj-truncation-iso} to embed any weighted projective space (or variety) into $\pee^N$ for some large enough $N$.
        As we've already mentioned, really this is a specific case of \cref{ssub:well_formed_weighted_projective_spaces}, since the weight $(1,\ldots,1)$ is well-formed, but in a sense it is also very different, since we are aiming for a \emph{specific} well-formed weighted projective space.
        We start this section with a technical lemma.

        \begin{lemma}\label{lem:suffic-large-d}
            Let\footnote{%
                This theorem is actually true for a general graded ring $R$, but we state it here in the more specific case, since it is the only one that we use.
            }
            \[
                R=\frac{\kazn}{I}
            \]
            for some radical weighted-homogeneous relevant ideal $I\triangleleft\kazn$.
            Then there exists some $d\in\nn$ such that $R^{(d)}$ is generated by $R_d$.
        \end{lemma}

        \begin{proof}
            We can simply apply the more general proof of \cite[Lemma~5.5.3]{Tevelev:jdT35_ao}.
        \end{proof}

        It turns out that \cref{lem:suffic-large-d} is the only new thing that we need to show that embedding into straight projective space is always possible.
        
        \begin{theorem}[Straightening out weighted projective space]\label{thm:straightening-wps}
            Let $X=\van(I)\subseteq\pazn$ be an irreducible\footnote{%
                We include the hypothesis that $X$ is irreducible here just to make the proof easier.
                In general, given some general weighted projective variety $V$, we can simply look at its irreducible components $V_i$ separately, which all embed into $\pee^{N_i}$ for some $N_i$.
                Then, since there is a natural embedding $\pee^N\hookrightarrow\pee^M$ for any $N\leqslant M$, we can embed $V$ into $\pee^N$ where $N=\max_i\{N_i\}$ by simply embedding $\pee^{N_i}\hookrightarrow\pee^N$ for each $i$.

                This argument does lack rigour at the end, when we simply `put all the pieces back together', but we do not have time to cover it here unfortunately.
                The author does not know of a suitable reference for this proof, but is sure that one must exist somewhere, for what that's worth.
            } non-empty weighted projective variety.
            Then there exists some $N$ large enough, and some projective variety $Y\subseteq\pee^N$, such that
            \[
                X\cong Y\subseteq\pee^N.\qedhere
            \]
        \end{theorem}

        \begin{proof}
            Let $R=\kazn/I$.
            Amongst other things, \cref{thm:proj-truncation-iso} tells us that, for any graded ring $S$ and any $d\in\nn$, prime weighted-homogeneous ideals in $S^{(d)}$ not containing $S^{(d)}_+$ are in exact correspondence with prime weighted-homogeneous ideals in $S$ not containing $S_+$.
            Thus, since \mbox{$I\triangleleft\kazn$} is a prime weighted-homogeneous relevant ideal ($X$ is non-empty) we know that it corresponds exactly to a prime weighted-homogeneous ideal $J\triangleleft\kazn^{(d)}$ not containing $\kazn^{(d)}_+$.

            That is,
            \[
                R = \frac{\kazn}{I} \implies R^{(d)} = \frac{\kazn^{(d)}}{J}
            \]
            where $J$ is prime weighted-homogeneous and such that $\van(J)\subseteq\proj\kazn^{(d)}$ is non-empty (since $J$ doesn't contain $\kazn^{(d)}_+$).\footnote{%
                This bit of the argument is admittedly a bit too hand-wavey for the author's liking.
                Hopefully looking at the example in \cref{sub:a_worked_example} will convince the reader that this could indeed be made more rigorous.
            }

            By \cref{lem:suffic-large-d} we can find some $d\in\nn$ such that $R^{(d)}$ is generated by $R_d$.
            In $R^{(d)}$ the elements of $R_d$ have degree $1$ by definition, and so
            \[
                R^{(d)} = \frac{k[y_0,\ldots,y_N]}{J}
            \]
            where $\wt y_i=1$ for all $i$, and $J$ is the image of $I$ in $R^{(d)}$, which is homogeneous (in the straight sense since $a=(1,\ldots,1)$) by the above argument.

            So \cref{thm:proj-truncation-iso} and \cref{thm:what-proj-means-for-quotients} tell us that
            \[
                X = \proj R \cong \proj R^{(d)} = \van(J)\subseteq\pee(1,\ldots,1) = \pee^N.\qedhere
            \]
        \end{proof}

        We showed in \cref{sub:morphisms_between_varieties} that our definition of isomorphisms agrees with the usual definition of isomorphisms between straight projective varieties.
        Thus we get the following corollary, which is really just \cref{thm:straightening-wps} phrased in a different way.

        \begin{corollary}\label{cor:wpv-can-be-though-of-as-spv-for-some-N}
            Let $X\subseteq\pazn$ be a non-empty weighted projective variety.
            Then $X$ can also be thought of as a straight projective variety inside some $\pee^N$.
        \end{corollary}

        


    \subsection{A worked example} 
    \label{sub:a_worked_example}
    
    To check that we have a working understanding of \cref{sub:truncation_of_graded_rings} we now look at an explicit example.
    This also gives us a chance to see how ideals transform under truncation, and to maybe help clarify the proof of \cref{thm:straightening-wps}.

    The example that we choose is from \cite[Exercise~5,~Section~6.5]{Tevelev:jdT35_ao}, and is a variation on \cite[Example~3.7]{Reid:2002uy}.

    Once again, we point out that we are not really covering the whole picture here -- we are dealing simply with the underlying topological spaces.
    As noted in \cref{thm:proj-truncation-iso}, the best source for the gory details is probably \cite[Proposition~(2.4.7)]{Grothendieck:1961tr}.

    \begin{example}\label{ex:that-worked-example}
        Let $f=x^5+y^3+z^2\in\cc_{(12,20,30)}[x,y,z]$, so that $f$ is weighted-homogeneous of degree $60$.
        Compute
        \[
            \proj\frac{\cc_{(12,20,30)}[x,y,z]}{(f)}.\qedhere
        \]
    \end{example}

    \begin{proof}[Solution]
        Since $f$ is weighted-homogeneous, $(f)$ is a weighted-homogeneous relevant ideal.
        Further, since $f$ is irreducible and $\cc_{(12,20,30)}[x,y,z]$ is a UFD\footnote{%
            If $R$ is a UFD then $R[t]$ is also a UFD.
        } the ideal $(f)$ is prime.
        By \cref{thm:what-proj-means-for-quotients} we know that this is the variety
        \[
            \van(x^5+y^3+z^2)\subseteq\pee(12,20,30),
        \]
        but \cref{thm:straightening-wps} makes us wonder what this looks like as a straight projective variety.

        Using the ideas in the proof of \cref{thm:can-always-well-form} we see that
        \[
            \underbrace{\pee(12,20,30)}_{\proj R} \cong \underbrace{\pee(6,10,15)}_{\proj R^{(2)}} \cong \underbrace{\pee(6,2,3)}_{\proj R^{(10)}} \cong \underbrace{\pee(3,1,3)}_{\proj R^{(20)}} \cong \underbrace{\pee(1,1,1)}_{\proj R^{(60)}}=\pee^2
        \]
        is a straightening of $\pee(12,20,30)=\proj\cc_{(12,20,30)}[x,y,z]$.
        So now we have to think how the ideal $(f)\triangleleft R$ transforms under these truncations.

        The isomorphism $\proj R\cong\proj R^{(d)}$ is induced by the inclusion $R^{(d)}\hookrightarrow R$ (using the ideas from \cref{sub:morphisms_between_varieties}).
        How canonical this isomorphism is, however, depends on whether we are looking at case (i) or case (ii) from \cref{thm:can-always-well-form}.
        That is, if $d\mid a_i$ for all $i$ then $R^{(d)}$ is simply $R$ with all the gradings divided through by $d$.
        So our weighted-homogeneous ideal $I\triangleleft\kazn$
        If, however, we have that $d\mid a_i$ for all $i\neq j$, and $\gcd(d,a_j)=1$, then the isomorphism is a tad less simple.
        It comes from the correspondence of prime ideals mentioned in our proof of \cref{thm:proj-truncation-iso}: weighted-homogeneous prime ideals $\prid\triangleleft R^{(d)}$ correspond uniquely to weighted-homogeneous prime ideals $\prid'\triangleleft R^{(d)}$ in such a way that, if $f\in \prid$, then $f^d\in \prid'$.

        So our first isomorphism is a very natural one, since $12,20,30$ all divide by $2$:
        \[
            \proj\frac{\cc_{(12,20,30)}[x,y,z]}{(f)} \cong \proj\left(\frac{\cc_{(12,20,30)}[x,y,z]}{(f)}\right)^{(2)} = \proj\frac{\cc_{(6,10,15)}[x,y,z]}{(f)}.
        \]
        But from now on, every isomorphism falls into case (ii) -- there is no common factor of \emph{all} of the $a_i$, only for pairs $a_i,a_j$.
        Let's look at the first such one.

        \bigskip

        By definition, $R^{(ab)}=(R^{(a)})^{(b)}$, and so we are interested first in
        \[
            \proj\bigg(\underbrace{\frac{\cc_{(6,10,15)}[x,y,z]}{(f)}}_{\text{call this } S}\bigg)^{(5)}.
        \]
        We first look at the `numerator' of this quotient.
        Since $y,z\in R_{5k}$ for all $k\in\nn$, and $\gcd(5,6)=1$, we see that
        \[
            \cc_{(6,10,15)}[x,y,z]^{(5)} = \cc_{(6,2,3)}[x^5,y,z].
        \]
        So let's turn now to the `denominator'.

        The ideal in $S^{(5)}$ corresponding to $(f)$ should be $(f^5)$, thus
        \[
            S^{(5)} = \frac{\cc_{(6,2,3)}[x^5,y,z]}{(f^5)},
        \]
        where $f$ is weighted-homogeneous of degree $30$ and so $f^5$ is weighted-homogeneous of degree $150$.
        But $\van(f)=\van(f^5)\subseteq\pee(6,10,15)$ by the more general fact that $\van(g)=\van(g^k)$ whenever\footnote{%
            We might not need such a strong condition on $g$, but we lose nothing here by erring on the side of caution.
        } $g$ is irreducible and $k\in\nn$.
        So \cref{thm:what-proj-means-for-quotients} tells us that
        \[
            \proj S^{(5)} = \proj\frac{\cc_{(6,2,3)}[x^5,y,z]}{(f^5)} = \proj\underbrace{\frac{\cc_{(6,2,3)}[x^5,y,z]}{(f)}}_{\text{call this }\overline{S^{(5)}}},
        \]
        which \emph{does} makes sense, as $f$ is weighted-homogeneous of degree $30$, and thus $f\in\cc_{(6,2,3)}[x^5,y,z]$ is still weighted-homogeneous, now of degree $6$.
        Finally for this first isomorphism, we simplify things a bit by using the isomorphism of graded rings
        \[
            \frac{\cc_{(6,2,3)}[x^5,y,z]}{(x^5+y^3+z^2)} \cong \frac{\cc_{(6,2,3)}[r,s,t]}{(r+s^3+t^2)}.
        \]

        Putting this all together gives us this composition of maps of rings that all have isomorphic $\proj$:
        \[
            \underbrace{\frac{\cc_{(6,10,15)}[x,y,z]}{(x^5+y^3+z^2)}}_{S} \mapsto \underbrace{\left(\frac{\cc_{(6,10,15)}[x,y,z]}{(x^5+y^3+z^2)}\right)^{(5)} = \frac{\cc_{(6,2,3)}[x^5,y,z]}{(x^5+y^3+z^2)^5}}_{S^{(5)}} \mapsto \underbrace{\frac{\cc_{(6,2,3)}[x^5,y,z]}{(x^5+y^3+z^2)} \cong \frac{\cc_{(6,2,3)}[r,s,t]}{(r+s^3+t^2)}}_{\overline{S^{(5)}}}.
        \]

        \bigskip

        Using the same process and notation as above, if we let $T=\overline{S^{(5)}}$ then we can repeat this with $T \mapsto T^{(2)} \mapsto \overline{T^{(2)}}$, where
        \[
            \overline{T^{(2)}} = \frac{\cc_{(3,1,3)}[u,v,w]}{(u+v^3+w)}.
        \]
        Finally, with $U=\overline{T^{(2)}}$ we have $U \mapsto U^{(3)} \mapsto \overline{U^{(3)}}$, where
        \[
            \overline{U^{(3)}} = \frac{\cc_{(1,1,1)}[X,Y,Z]}{(X+Y+Z)}.
        \]
        So, at long last, we see that
        \[
            \proj\frac{\cc_{(12,20,30)}[x,y,z]}{(x^5+y^3+z^2)} \quad\cong\quad \proj\frac{\cc_{(1,1,1)}[X,Y,Z]}{(X+Y+Z)} \quad=~ \underbrace{\van(X+Y+Z)}_{\text{sitting inside }\pee(1,1,1)=\pee^2} ~\cong\quad \pee^1.\qedhere
        \]
    \end{proof}

    Now, having done this example, let's talk through what it means.
    First of all, although we have shown that our original weighted projective variety $X=\van(x^5+y^3+z^2)\subset\pee(12,20,30)$ is isomorphic to $\pee^1$, it doesn't mean that it's exactly the same.
    The affine cone over $\pee^1$ is simply the plane $\aff^2$, whereas the affine cone $\hat{X}$ over $X$ is a degree-$5$ hypersurface inside $\aff^3$.
    It turns out in fact that $X$ is a rather special singularity, see \cite[Example~3.7]{Reid:2002uy} for more information, since it isn't too relevant here (but it is very interesting).\footnote{%
        See also the answer to one of the author's questions (which also displays their original misunderstandings) for a discussion about the links to the Poincaré Homology Sphere: \cite{1426420}.
    }

    Another thing to note is that we travelled down the algebraic path in our solution of \cref{ex:that-worked-example}, but as we might have expected, we could have instead followed a more geometric one.
    \cref{eq:affine-cone-thing} tells us that $X$ could be thought of as the quotient variety
    \[
        X = \frac{\vanaff(x^5+y^3+z^2)}{\gm} \subseteq \frac{\aff^3\noz}{\gm},
    \]
    and so we might have tried to construct this quotient explicitly in an attempt to understand the structure of $X$.
    We chose not to, because here the algebraic approach gives us a nice way of using all the things that we've found out so far, and because it is arguably much slicker.

    But we could also think of weighted projective space as a quotient of straight projective space, and construct $X$ as a quotient of a straight projective variety.
    What do we mean by this?
    Well, we claim that we have quotient maps
    \begin{equation}\label{eq:quotient-diagram}
        \begin{tikzcd}[column sep=huge, row sep=6em]
            & \pee^n \dar[two heads]{\bigoplus_{i=0}^n\mu^{a_i}} \\
            \aff^{n+1}\noz \rar[two heads]{\gm^{(a)}} \urar[two heads]{\gm^{(1,\ldots,1)}} & \pazn
        \end{tikzcd}
    \end{equation}
    where the labels on the arrows are the groups that we quotient by to obtain the surjection.
    We study at a slightly more specific case of this in \cref{sec:curves_in_weighted_projective_space} (looking only at plane curves) but the general idea stays the same.
    Up until now we have been studying the quotient along the bottom of \cref{eq:quotient-diagram}, and the upper-left quotient is a specific case of it, giving the usual well-understood case of projective space.
    The quotient on the right hasn't really been mentioned at all yet, but it comes in useful if we want to use what facts we know about straight projective space (which are sometimes `nicer' than those about affine space) to gain some know-how about weighted projective spaces.



\section{Plane curves in weighted projective space} 
\label{sec:curves_in_weighted_projective_space}

    {\color{purple}\textbf{The following section has been edited due to the discovery of a mistake in the proofs of \cref{lem:cover-is-smooth-pc-and-things,lem:that-smoothness-lemma}, first noticed by Claude Quitte, and worked out in detail by Bernhard Albach. We have left the offending lemmas in, so as not to nullify any existing references to them, but have changed their contents to instead contain certain \emph{hypotheses}. For a more general solution to the problem (by a different construction of a straight cover that actually \emph{is} non-singular), see the upcoming thesis of Bernhard Albach.}}

    \bigskip

    \emph{In this section we assume familiarity with some of the fundamentals of Riemann surfaces and maps between them (see \cite[Chapters~1,~2]{Miranda:1995uz}).}
    \emph{We also assume knowledge of orbit spaces, covering spaces, and some other concepts from topology, though these are usually explained when used.}
    \emph{Finally, we assume some facts about algebraic curves, though these are stated before being used.}

    \bigskip

    \begin{note}
        In this section (\cref{sec:curves_in_weighted_projective_space}) \emph{only} we write $\pee=\pathree$, $\poly=\kathree$, and $f(x)=f(x_0,x_1,x_2)$.
        But the weight \mbox{$a=(a_0,a_1,a_2)$} and the indeterminates $x_0,x_1,x_2$ are still lurking about in the background, and when we say `weighted-homogeneous' we still mean `$a$-weighted-homogeneous'.
        Note that when we write $\pee^k$ we mean projective $k$-space, as per usual, and we will write $\pee^1$ for the projective line, so there should be no ambiguity in our use of $\pee$ to denote $\pee(a_0,a_1,a_2)$.
    \end{note}
    Now seems like a good time to take a look at a specific family of weighted projective varieties, namely \emph{plane curves}, so that we don't lose our geometric intuition, and so that we have some (comparatively) concrete examples to hold on to and examine.
    As mentioned at the end of \cref{sec:algebra}, we will half-break our promise made in \cref{sub:notation_and_conventions}: we assume some prior knowledge of straight projective algebraic geometry, and so don't deal with it as a special case of weighted projective algebraic geometry.
    Instead, we now shift our viewpoint slightly to think of weighted projective space as a quotient of straight projective space.
    By doing so, we also get a change to take a break from the algebra side of things and to work primarily with the geometry of these objects that we're studying.

    An important thing to be aware of is that the path we follow here is almost certainly much much longer than it needs to be, but (in the view of the author) we uncovers plenty of nice facts along the way, and also develop multiple ways of viewing these objects.

    \begin{definition}[Plane curves in weighted projective space]\label{defn:plane-curves-in-wps}
        Let $f=f(x_0,x_1,x_2)\in\poly$ be a weighted-homogeneous degree $d$ polynomial with no repeated factors.\footnote{%
            The reason that we say this is really just to simplify things without losing generality.
            We already know that $\van(f^k)=\van(f)$, and in a similar way we see that $\van(f^k g)=\van(fg)$.
            So we might as well assume that our polynomial has no repeated factors, but it isn't entirely necessary for our purposes.
        }
        Then
        \[
            C_f = \van(f) \subseteq \pee
        \]
        is a \emph{degree-$d$ plane curve in $\pathree$}.

        We say that a plane curve $C$ is \emph{irreducible} if $f$ has no non-constant factors apart from scalar multiples of itself, since then $C$ cannot be written as a non-trivial union of other plane curves (using \cref{lem:topology-well-defined}).
    \end{definition}

    \begin{definition}[Singular points]
        Let $f\in\kathree$ be a degree-$d$ weighted-homogeneous polynomial.
        Then we say that $p=(p_0,p_1,p_2)\in\aff^3\noz$ is a \emph{singular point of $f$} if
        \[
        \eval{\pdv{f}{x_0}}_p = \eval{\pdv{f}{x_1}}_p = \eval{\pdv{f}{x_2}}_p = 0.
        \]
        We say that $f$ is \emph{non-singular} if it has no singular points.
        Similarly, we say that a plane curve $C=C_f$ is \emph{non-singular} if its defining polynomial\footnote{%
            So if we dropped the requirement that polynomials have no repeated factors then we'd need to specify \emph{which} of the infinitely-many defining polynomials we mean -- namely the one with no repeated factors.
        } $f$ is either non-singular, or singular only at points outside of $C$, i.e. only at points $p$ such that $f(p)\neq0$.
    \end{definition}

    When $a=(1,1,1)$, i.e. in the straight case, we see that our definition of plane curves is exactly the same as the usual definition for projective plane curves.
    We now state a fundamental fact about plane curves in straight projective space.

    \begin{lemma}\label{lem:str-proj-plane-curve-rs}
        Let $C\subset\pee^2$ be a non-singular plane curve.
        Then $C$ is a compact Riemann surface.
    \end{lemma}

    \begin{proof}
        See \cite[Chapter~I,~Proposition~3.6]{Miranda:1995uz}.
    \end{proof}

    \subsection{Some facts about different notions of quotients} 
    \label{sub:some_facts_that_we_need}

        In this subsection we state, and cite proofs for, a lot of technical lemmas that we will use in \cref{sub:weighted_projective_plane_curves_as_riemann_surfaces}.
        We can split these into three types by looking at what notion of `quotient' they concern: topological spaces, Riemann surfaces, and projective GIT quotients.

        We have a few notational notes before we start that seemed rather pointless to put in \cref{ssub:list_of_notation_and_assumptions} since they are only really used here.
        Given a group action $G$ on some space $X$ we write $G_p$ for a point $p\in X$ to mean the stabiliser subgroup $\{g\in G\mid g\cdot p = p\}$.
        We write $\sigma_n$ to mean a general $n$-th root of unity (so $\sigma_n=\omega_n^k$ for some $0\leqslant k<n$).
        Finally, with a holomorphic map $f\colon X\to Y$ between Riemann surfaces we write $\mult_p(f)$ to mean the ramification index, as in \cite[Chapter~II,~Definition~4.2]{Miranda:1995uz}.

        \subsubsection{Topological quotients (orbit spaces)} 
        \label{ssub:topological_quotients_}

            \begin{lemma}[Quotients preserve compactness]\label{lem:comp-conn-quot}
                Let $X$ be a compact topological space and $G$ some finite group acting on $X$.
                Then the orbit space $X/G$ is a compact topological space.
            \end{lemma}

            \begin{proof}
                This follows from the standard fact that a continuous image of a compact space is compact, and that the quotient map is continuous (by definition of the quotient topology on the quotient space).
            \end{proof}
        

        \subsubsection{Quotients of Riemann surfaces by group actions} 
        \label{ssub:quotients_of_riemann_surfaces_by_group_actions}

            \begin{definition}[Group actions on a Riemann surface]
                Let $G$ be a group acting on a Riemann surface $X$.
                Then we say that the action of $G$ is
                \begin{itemize}
                    \item \emph{holomorphic} if the bijection $\varphi_g\colon X\to X$ given by $x\mapsto g\cdot x$ is holomorphic for all $g\in G$;
                    \item \emph{effective} if the kernel $K=\{g\in G \mid g\cdot x = x\text{ for all }x\in X\}$ is trivial.\qedhere
                \end{itemize}
            \end{definition}

            \begin{theorem}[Quotient of a Riemann surface by a group action]\label{thm:group-on-rs}
                Let $G$ be a finite group acting holomorphically and effectively on a Riemann surface $X$.
                Then we can endow the orbit space $X/G$ with the structure of a Riemann surface.
                Moreover, the quotient map $\vartheta\colon X\to X/G$ is holomorphic of degree $|G|$ and $\mult_x(\vartheta)=|G_x|$ for any point $x\in X$.
            \end{theorem}

            \begin{proof}
                See \cite[Chapter~III,~Theorem~3.4]{Miranda:1995uz}.
            \end{proof}
        

        \subsubsection{Projective GIT quotients} 
        \label{ssub:projective_git_quotients}

            All of the following is taken from \cite{Hoskins:2012uq}, but has been phrased here slightly differently, and in a different order, just to avoid getting too carried away with the vast subject of geometric invariant theory (GIT).
            We also oversimplify the background machinery wildly, so do check \cite[Chapter~4]{Hoskins:2012uq} for the whole story.

            \begin{definition}[Linear action of reductive groups on projective varieties]
                Let a reductive\footnote{We are only ever interested in finite groups here though, and all finite groups are reductive.} group $G$ act on a projective variety $X\subseteq\pee^n$.
                Then its action is said to be \emph{linear} if $G$ acts via a homomorphism $G\to\GL(n+1)$.
            \end{definition}

            \begin{definition}[Projective GIT quotient {(\cite[Definition~4.6]{Hoskins:2012uq})}]\label{defn:proj-git-quotient}
                Let $G$ be a reductive group with a linear action on a projective variety $X\subset\pee^n$.
                We define the \emph{projective GIT quotient variety $X\sslash G$} to be the projective variety given by $\proj S(X)^G$, where $S(X)$ is the homogeneous coordinate ring of $X$.
            \end{definition}

            \begin{definition}[Geometric quotient {(\cite[Definition~2.28]{Hoskins:2012uq})}]\label{defn:geometric-quotient}
                A projective GIT quotient map\footnote{Really the map is only defined on $X^{\text{ss}}\subset X$, where $X^{\text{ss}}$ is a certain subset of $X$. Again though, this is all much beyond what is needed here.} $\phi\colon X\twoheadrightarrow X\sslash G$ is said to be \emph{geometric} if the preimage $\phi^{-1}([x])$ of each point $[x]\in X\sslash G$ is a single orbit in $X$.
            \end{definition}

            So if $\phi\colon X\twoheadrightarrow X\sslash G$ is a geometric projective GIT quotient then $X\sslash G$ is simply the topological quotient (i.e. the orbit space) $X/G$.
            In particular then, $X\sslash G$ naturally has the quotient topology coming from the quotient map $X\to X/G$.



    \subsection{Weighted projective plane curves as Riemann surfaces} 
    \label{sub:weighted_projective_plane_curves_as_riemann_surfaces}

        \emph{We now make the assumption that our weighted projective space $\pee$ is well-formed so that, in particular, the $a_i$ are pairwise coprime.}
        \emph{Note that we don't lose much generality by doing this, since \cref{thm:can-always-well-form} tells us that any weighted projective space is isomorphic to a well-formed one.}
        \emph{The only information that we lose by passing to this isomorphic copy is the specific embedding of our original variety, but here we are much less interested in the embedding of plane curves and much more so in their intrinsic nature.}

        \begin{lemma}\label{lem:pi-hash-maps-homo-to-homo}
            Define\footnote{%
                The reason for this choice of notation is that we will eventually show that $\pi_{\#}$ is the pushforward of a quotient map $\pi$.
            } the homomorphism of graded rings
            \begin{align*}
                \pi_{\#}\colon\kathree&\to k[y_0^{a_0},y_1^{a_1},y_2^{a_2}] \\
                x_i&\mapsto y_i^{a_i}.
            \end{align*}
            Let $f\in\kathree$ be weighted-homogeneous of degree $d$.
            Then $\pi_\#(f)$ is homogeneous of degree $d$.
        \end{lemma}

        \begin{proof}
            This follows from the definition of the degree of a weighted-homogeneous polynomial (\cref{defn:wh-poly}).
        \end{proof}

        \begin{note}
            We try to be consistent with the convention that the $a$-weighted polynomial ring has indeterminates $x_i$ and the usual polynomial ring has indeterminates $y_i$, but sometimes we slip up, and sometimes for a good reason (to avoid unnecessary complication with notation at times).
            It is just a choice of notation, so doesn't affect that maths at all, but can be confusing.
            Just be aware.
        \end{note}

        \begin{definition}[Straight cover]
            Given some plane curve $C=C_f\in\pee$ we define its \emph{straight cover $\cover{C}$} as the straight projective variety
            \[
                \cover{C} = C_{\cover{f}} = \van(\pi_\#(f))\subset\pee^2.
            \]
            Note that this map is well defined by \cref{lem:pi-hash-maps-homo-to-homo}.
        \end{definition}

        We need to explain why we give this variety $\cover{C}$ the name `straight cover', since it is a very suggestive name.
        The rest of this section is sort of dedicated to showing why we use this name.

        \begin{example}
            Let $f(x,y,z)=x^4+y^4+z^2+xyz\in k_{(1,1,2)}[x,y,z]$ be weighted-homogeneous of degree-$4$, giving us the plane curve
            \[
                C = C_f = \van(x^4+y^4+z^2+xyz) \subset \pee(1,1,2).
            \]
            Then $\cover{f} = \pi_\#(f) = x^4+y^4+z^4+xyz^2\in k[x,y,z]$ is also degree-$4$ homogeneous, and
            \[
                \cover{C} = C_{\cover{f}} = \van(x^4+y^4+z^4+xyz^2) \subset \pee^2.\qedhere
            \]
        \end{example}

        It turns out that, for a nice enough plane curve $C_f$, the straight cover $C_{\cover{f}}$ is non-singular -- a fact which we now state and prove, as well as saying what exactly we mean by `nice enough'.

        \begin{lemma}\label{lem:cover-is-smooth-pc-and-things}
            {\color{purple} The original content of this lemma was as follows:}
            \begin{quote}
            Let $f\in\kathree$ be weighted-homogeneous and non-singular.
            Then $\cover{f}=\pi_\#(f)$ is homogeneous and non-singular.
            Equivalently, if we have some non-singular plane curve $C_f\subset\pee$ then its straight cover $C_{\cover{f}}\subset\pee^2$ is also a non-singular plane curve.
            \end{quote}
            {\color{purple}This is not, however, in general, true. \textbf{Instead, from now on, we always assume the additional hypotheses that $C$ is such that $\cover{C}$ is non-singular.}}
        \end{lemma}

        \begin{lemma}\label{lem:quotient-of-cover}
            Let $C=C_f\subset\pee$ be a non-singular plane curve and $\cover{C}$ its straight cover.
            Let $G=\mu^{a_0}\times\mu^{a_1}\times\mu^{a_2}$ and define an action of $G$ on $\cover{C}$ by
            \[
                g\cdot y = (\sigma_{a_0},\sigma_{a_1},\sigma_{a_2})\cdot [y_0:y_1:y_2] = [\sigma_{a_0}y_0:\sigma_{a_1}y_1:\sigma_{a_2}y_2].
            \]
            Then
            \begin{enumerate}[(i)]
                \item the orbit space $\cover{C}/G$ is a compact Riemann surface;
                \item the quotient map $\vartheta\colon \cover{C}\twoheadrightarrow\cover{C}/G$ is holomorphic of degree $a_0a_1a_2$;
                \item $\mult_y(\vartheta)=|G_y|$ for any point $y\in\cover{C}$.\qedhere
            \end{enumerate}
        \end{lemma}

        \begin{proof}
            We have already done most of the hard work for this proof, so we just need to fit all the pieces together.

            \cref{lem:str-proj-plane-curve-rs,lem:cover-is-smooth-pc-and-things} say that $\cover{C}$ is a compact Riemann surface, so \cref{lem:comp-conn-quot} tells us that $Y/G$ is a compact topological space.
            \cref{thm:group-on-rs} gives us the rest of the claims \emph{assuming that we can show that} $G$ acts holomorphically and effectively (since it is finite of order $a_0a_1a_2$).

            The kernel of the action of $G$ on $\cover{C}$ is
            \[
                K = \{g\in G\mid g\cdot y = y\text{ for all }y\in \cover{C}\}
            \]
            but since we have assumed that $a$ is well-formed, and thus that the $a_i$ are all pairwise coprime, we know that $\mu^{a_i}\cap\mu^{a_j}=\{1\}$ for $i\neq j$.
            By definition, $g\cdot y=y$ for all $y\in\cover{C}$ if and only if $g=(\lambda,\lambda,\lambda)$ for some $\lambda\in k\noz$.
            So the above comment tells us that no $g\in G$ is of this form apart from $(1,1,1)$, which is the identity in $G$.
            Thus the action is effective, since the kernel is trivial.
            As for the map $\varphi_g$ being holomorphic, this follows straight away from the fact that it is an algebraic map.
            More specifically, it is simply a polynomial map with constant coefficients.
        \end{proof}

        \begin{lemma}\label{lem:git-quotient-exists}
            Let $C=C_f\subset\pee$ be a non-singular plane curve and $\cover{C}\subset\pee^2$ be its straight cover, which is also a non-singular plane curve by \cref{lem:cover-is-smooth-pc-and-things}.
            Let $G=\mu^{a_0}\times\mu^{a_1}\times\mu^{a_2}$ act on $\cover{C}$ as in \cref{lem:quotient-of-cover}.
            Then this induces an action of $G$ on the homogeneous coordinate ring $S(\cover{C})$, and we have a well-defined GIT quotient
            \[
                \varphi\colon\cover{C}=\proj S(\cover{C})\twoheadrightarrow\proj S(\cover{C})^G=\cover{C}\sslash G.\qedhere
            \]
        \end{lemma}

        \begin{proof}
            By \cref{defn:proj-git-quotient}, we need to show that $G$ is reductive and has linear action on the homogeneous coordinate ring $S(\cover{C})$.
            Now $G$ is reductive by definition, since it is finite.
            Further, the action is linear, since it acts diagonally.
            That is,
            \[
                (\sigma_{a_0},\sigma_{a_1},\sigma_{a_2})\cdot[y_0:y_1:y_2] = \left[\mqty(\dmat[0]{\sigma_{a_0},\sigma_{a_1},\sigma_{a_2}})\mqty(y_0\\y_1\\y_2)\right].
            \]
            So we have the well-defined GIT quotient $\cover{C}\sslash G = \proj S(\cover{C})^G$, but what does this look like?

            The induced action of $G$ on $k[y_0,y_1,y_2]$ is given by
            \[
                (\sigma_{a_0},\sigma_{a_1},\sigma_{a_2})\cdot f(y_0,y_1,y_2) = f(\sigma_{a_0}y_0,\sigma_{a_1}y_1,\sigma_{a_2}y_2).
            \]
            Recall that we have assumed that $a$ is well formed, and hence that the $a_i$ are all pairwise coprime, so $\mu^{a_i}\cap\mu^{a_j}=\{1\}$ for $i\neq j$.
            Thus
            \[
                S(\cover{C})^G = \left(\frac{k[y_0,y_1,y_2]}{(\cover{f})}\right)^G = \frac{k[y_0^{a_0},y_1^{a_1},y_2^{a_2}]}{(\cover{f})}
            \]
            since $\cover{f}$ is already a polynomial in $y_i^{a_i}$ by definition, and
            \[
                \cover{C}\sslash G = \proj \frac{k[y_0^{a_0},y_1^{a_1},y_2^{a_2}]}{(\cover{f})}.\qedhere
            \]
        \end{proof}

        This section so far has been not much more than a wall of text consisting solely of definitions, lemmas, and proofs, so let's now take a break and have a look at what we've actually discovered and defined.

        Given some non-singular plane curve $C=C_f\subset\pee$, by using its straight cover $\cover{C}=C_{\cover{f}}$ and the finite group $G=\mu^{a_0}\times\mu^{a_1}\times\mu^{a_2}$ we can build two more objects, giving us three in total, including $C$.
        We later claim (\cref{cor:quotient-set-up-with-pi-etc}) that all three actually give us the same thing, and so we have three ways of looking at non-singular plane curves.
        The three objects we have are
        \begin{enumerate}[(i)]
            \item the non-singular plane curve $C=C_f\subset\pee$;
            \item a quotient map $\vartheta\colon\cover{C}\twoheadrightarrow\cover{C}/G$ of compact Riemann surfaces (\cref{lem:quotient-of-cover}));
            \item a GIT quotient map $\varphi\colon\cover{C}\twoheadrightarrow\cover{C}\sslash G$ of projective varieties (\cref{lem:git-quotient-exists}).
        \end{enumerate}

        \begin{theorem}\label{thm:i-and-iii-equiv}
            The objects (i) and (iii) are equivalent.
            That is, $C\cong\cover{C}\sslash G$ as weighted projective varieties.
        \end{theorem}

        \begin{proof}
            It turns out that, not only are these varieties isomorphic, they are `nicely' isomorphic.
            That is, the isomorphism of varieties arises from an isomorphism of graded rings, namely
            \begin{align*}
                S(C) = \frac{\kathree}{(f)} &\rightarrow \frac{k[y_0^{a_0},y_1^{a_1},y_2^{a_2}]}{(\cover{f})} = S(\cover{C}\sslash G)\\
                x_i &\mapsto y_i^{a_i}.
            \end{align*}
            This induces the desired isomorphism of weighted projective varieties $C\cong\cover{C}\sslash G$.
            \emph{In this section we denote this isomorphism by $\psi\colon \cover{C}\sslash G\to C$, where
            \[
                \psi\colon\orbit_G([p_0:p_1:p_2])\mapsto[p_0^{a_0}:p_1^{a_1}:p_2^{a_2}].\qedhere
            \]}
        \end{proof}

        \begin{lemma}\label{lem:git-quotient-is-orbit-space}
            The objects (iii) and (ii) are equivalent.
            That is, the GIT quotient $\cover{C}\sslash G$ is exactly the orbit space $\cover{C}/G$, which has all the structure of a compact Riemann surface.
        \end{lemma}

        \begin{proof}
            All that this lemma is really saying is that the quotient $\cover{C}\sslash G$ is geometric, as defined in \cref{defn:geometric-quotient}.
            So we need to show that the preimage of each point in $\cover{C}\sslash G$ is a single orbit in $\cover{C}$.

            Let $p\in\cover{C}\sslash G$ be a point.
            By \cref{thm:i-and-iii-equiv} we know that we can think of $\cover{C}\sslash G$ as the plane curve $C\subset\pee$ using the isomorphism $\psi\colon\cover{C}\sslash G\to C$, and so we can think of $p$ as the point $\psi(p)=|p_0:p_1:p_2|\in\pee$.
            Then
            \begin{align*}
                \varphi^{-1}(p) &= (\varphi^{-1}\circ\psi^{-1})(|p_0:p_1:p_2|) \\
                &= \{[\sigma_0 p_0^{1/a_0}:\sigma_1 p_1^{1/a_1}:\sigma_2 p_2^{1/a_2}] \mid \sigma_i\in\mu^{a_i}\} \\
                & = \{g\cdot[p_0^{1/a_0}:p_1^{1/a_1}:p_2^{1/a_2}] \mid g\in G\}
            \end{align*}
            is the orbit of a single point in $\cover{C}$, namely $[p_0^{1/a_0}:p_1^{1/a_1}:p_2^{1/a_2}]$.
        \end{proof}

        We can now state and prove the main result of \cref{sub:weighted_projective_plane_curves_as_riemann_surfaces}, which happens to be no more than a corollary of all the technical heavy lifting we've done already.

        \begin{corollary}\label{cor:quotient-set-up-with-pi-etc}
            Let $C=C_f\subset\pee$ be a non-singular plane curve and $\cover{C}\subset\pee^2$ its straight cover.
            Then the map $\pi\colon\cover{C}\to C$ given by
            \[
                \pi\colon[y_0:y_1:y_2]\mapsto|y_0^{a_0}:y_1^{a_1}:y_2^{a_2}|
            \]
            is a surjective map of compact Riemann surfaces.
            Further, $\pi$ is holomorphic of degree $a_0a_1a_2$ and such that $\mult_y(\pi)=|G_y|$ for any $y\in\cover{C}$.
        \end{corollary}

        \begin{proof}
            \emph{Although there is a lot of notation here, due to all these isomorphisms and quotient maps, the idea behind this is very simple, and has been pretty much explained by what we have done so far.
                        This is just putting all the pieces together and chasing notation around.}

            First we look at the map $\pi_\#\colon\kazn\to k[y_0^{a_0},y_1^{a_1},y_2^{a_2}]$ from \cref{lem:pi-hash-maps-homo-to-homo}.
            We can represent $\pi_\#$ by polynomials $\Pi_i\in k[y_0^{a_0},y_1^{a_1},y_2^{a_2}]$ where $\Pi_i=\pi(x_i)=y_i^{a_i}$.
            This induces a map $\pi\colon\proj k[y_0^{a_0},y_1^{a_1},y_2^{a_2}]\to\proj\kazn$ given by
            \[
                \pi\colon p=|p_0:p_1:p_2|\mapsto|\Pi_0(p):\Pi_1(p):\Pi_2(p)| = |p_0^{a_0}:p_1^{a_1}:p_2^{a_2}|.
            \]

            \cref{lem:git-quotient-is-orbit-space} says that $\cover{C}\sslash G=\cover{C}/G$, and so $\varphi=\vartheta$, since they both map a point in $\cover{C}$ to its $G$-orbit.
            Then we use the isomorphism $\psi\colon\cover{C}/G=\cover{C}\sslash G\twoheadrightarrow C$ to define the composition $\vartheta\circ\psi$ which has all the properties of $\vartheta$ from \cref{lem:quotient-of-cover}.
            But
            \[
                \vartheta\circ\psi\colon [y_0:y_1:y_2]\mapsto\orbit_G([y_0:y_1:y_2])\mapsto|y_0^{a_0}:y_1^{a_1}:y_2^{a_2}|.
            \]
            Thus $\vartheta\circ\psi=\pi$, and using \cref{lem:quotient-of-cover}, $\pi$ is a holomorphic map of degree $a_0a_1a_2$ between compact Riemann surfaces such that $\mult_y(\pi)=|G_y|$ for any $y\in\cover{C}$.
        \end{proof}

        Figure~\ref{eq:all-the-maps} is intended to be an understandable summary of all the confusing notation being thrown around, and the proof of \cref{cor:quotient-set-up-with-pi-etc} essentially aims to prove that the diagram commutes.

        
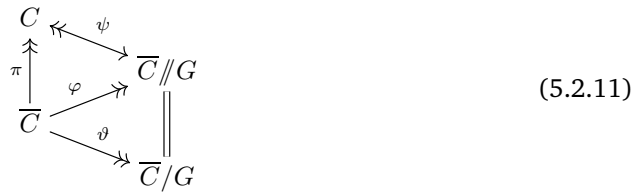
\begin{equation}\label{eq:all-the-maps}
            \begin{tikzcd}[column sep=3em, row sep=0.5em]
                C
                &
                \\
                &
                \cover{C}\sslash G
                    \ular[hookleftarrow, two heads, swap]{\psi}
                    \arrow[-, double equal sign distance]{dd}
                \\
                \cover{C}
                    \arrow[two heads]{uu}{\pi}
                    \urar[two heads]{\varphi}
                    \drar[two heads]{\vartheta}
                &
                \\
                &
                \cover{C}/G
            \end{tikzcd}
        \end{equation}

        From all of the above we also get a bonus corollary for free.
        Even though \cref{cor:wpv-can-be-though-of-as-spv-for-some-N} can give us the same result, we mention it here anyway, just to show that this could be an alternative path of getting to it.

        \begin{corollary}
            Let $C\subset\pee$ be a non-singular plane curve.
            Then $C$ is (isomorphic to) a projective variety.
        \end{corollary}

        \begin{proof}
            We know that $C\cong\cover{C}\sslash G$ by \cref{thm:i-and-iii-equiv}, and projective GIT quotients are, in particular, projective varieties (\cref{defn:proj-git-quotient}).
        \end{proof}

        \subsubsection{What happens to the affine patches?} 
        \label{ssub:what_happens_to_the_affine_patches_}

        We briefly discuss the story of the affine patches now, though we don't dedicate too much time to it, since it is slightly irrelevant compared to the results we move on to state and prove for the rest of the paper.
        But it is still interesting enough to be worth a mention.

        It is natural to think that there should be some way of `ungluing' a plane curve $C\subset\pee$ into its three affine patches, and then gluing them together in some standard way (homogenising the polynomials in such a way that they agree) to obtain the straight cover $\cover{C}$.
        This is slightly complicated, though, by the fact that the covering affine patches map onto the quotient affine patches by $\pi_i$ (that is, the quotient map for a $\mu^{a_i}$ action), but the straight cover maps onto the plane curve by $\pi=\pi_{ijk}$ (that is, the quotient map for a $\mu^{a_0}\times\mu^{a_1}\times\mu^{a_2}$ action).
        In a sense, we have to split $C$ up into its quotient affine patches, map back up to the covering affine patches, \emph{and then factor through some sort of patches, also in the affine plane}, before finally gluing them back together.

        Writing $D_i$ to mean the quotient affine patches, $\othercover{D_i}$ to mean the covering affine patches, $\cover{D}_i$ to mean the covers of these affine patches we can draw a diagram of the situation: see Figure~\ref{eq:notation-diagram}.
        We hope that the diagram is at least reasonably helpful and understandable, though we do stress its irrelevance to what is to follow.

        \begin{equation}\label{eq:notation-diagram}
            \begin{tikzcd}[column sep=1em, row sep=1.2em]
                &
                \textcolor{lightgray}{\aff^2}
                    \arrow[two heads, near start, lightgray]{dd}{\pi_{jk}}
                    \arrow[hookrightarrow, two heads, lightgray, near start]{rr}{\sim}
                &
                &
                \textcolor{lightgray}{V_i}
                    \arrow[hookrightarrow, near start, lightgray]{rr}{\kappa}
                &
                &
                \textcolor{lightgray}{\pee^2}
                    \arrow[two heads, near start, lightgray]{dddd}{\pi}
                &
                &
                \\
                \cover{D}_i
                    \arrow[-, darkgray]{ur}
                    \arrow[two heads, near start,]{dd}{\pi_{jk}}
                    \arrow[hookrightarrow, two heads, near start]{rr}{\sim}
                &
                &
                \cover{C}_i
                    \arrow[-, darkgray]{ur}
                    \arrow[hookrightarrow, near start]{rr}{\kappa}
                &
                &
                \cover{C}
                    \arrow[-, darkgray]{ur}
                    \arrow[two heads, near start]{dddd}{\pi}
                &
                &
                &
                \\
                &
                \textcolor{lightgray}{\aff^2}
                    \arrow[two heads, near start, lightgray]{dd}{\pi_i}
                &
                &
                &
                &
                &
                &
                \\
                \othercover{D}_i
                    \arrow[-, darkgray]{ur}
                    \arrow[two heads, near start,]{dd}{\pi_i}
                &
                &
                &
                &
                &
                &
                &
                \\
                &
                \textcolor{lightgray}{A_i}
                    \arrow[hookrightarrow, two heads, lightgray, near start]{rr}{\sim}
                &
                &
                \textcolor{lightgray}{U_i}
                    \arrow[hookrightarrow, near start, lightgray]{rr}{\iota}
                &
                &
                \textcolor{lightgray}{\pee}
                &
                &
                \textcolor{lightgray}{\bullet}
                \\
                D_i \arrow[-, darkgray]{ur}
                    \arrow[hookrightarrow, two heads, near start]{rr}{\sim}
                &
                &
                C_i \arrow[-, darkgray]{ur}
                    \arrow[hookrightarrow, near start]{rr}{\iota}
                &
                &
                C
                    \arrow[-, darkgray]{ur}
                &
                &
                \bullet
                    \arrow[hookrightarrow, xshift=-1ex, darkgray]{ur}{\text{inclusion}}
                    \arrow[twoheadleftarrow, swap, xshift=1ex, darkgray]{ur}{\text{restriction}}
                &
            \end{tikzcd}
        \end{equation}



    \subsection{The degree-genus formula for weighted projective varieties} 
    \label{sub:degree_genus}

        \emph{In this subsection we often talk of whether or not a polynomial has an $x_i$ term, or an $x_ix_j^k$ monomial, or something similar.
        When we say this, we implicitly mean a non-zero term, whether or not we mention it explicitly.
        So if, for example, we say that $f$ has a $x_i$ term, then we mean that $f=\lambda x_i+g$ where $\lambda\in k\noz$ and $g$ is some other polynomial in the $x_i$.}

        \bigskip

        \subsubsection{What does sufficiently general mean?} 
        \label{ssub:what_does_sufficiently_general}

            In the usual study of plane curves one tends to ignore the edge cases where a certain class of curve is poorly behaved.
            For example, all conics in $\pee^2$ are equivalent via a projective transformation to $x^2+y^2+z^2$, \emph{apart from} the singular cases that are equivalent to one of $x^2$ or $x^2+y^2$.
            In a sense, it's natural to think that these two examples will be difficult, because we are trying to study degree-$2$ polynomials in $k[x,y,z]$, but $x^2+y^2$ doesn't even have a non-zero $z$ term.
            It looks like it might be more at home as being classed as a degree-$2$ polynomial in $k[x,y]$, and $x^2$ is an even more extreme case.
            None of this is particularly rigorous, but it gives us a good idea of the restrictions we might want to place upon our polynomials to ensure that they define sufficiently nice plane curve.
            That is, we want our polynomials to be \emph{sufficiently general} in some sense.

            \begin{definition}[Sufficiently general]\label{defn:suff-general}
                A degree-$d$ weighted-homogeneous polynomial $f\in\kathree$ is \emph{sufficiently general} if $f$ satisfies the following for each $i$:
                \begin{enumerate}[(i)]
                    \item if $a_i\mid d$ then $f$ contains an $x_i^{d/a_i}$ term;
                    \item if $a_i\nmid d$ then $f$ contains an $x_jx_i^m$ term, where $j\neq i$ and $m=(d-a_j)/a_i$.
                \end{enumerate}
                We \emph{also} place some restrictions on the values of $d$ and all the $a_i$:
                \begin{itemize}
                    \item $d\geqslant2$;\footnote{%
                        Linear polynomials can only involve monomials $x_i$ where $a_i=1$, and so we just end up studying linear polynomials in straight projective space.
                        Hence we don't really lose too much generality in excluding these cases.
                    }
                    \item $d\geqslant a_i$;\footnote{%
                        So in the cases where $a_j\neq1$ for some $j$ (that is, all of the not-straight cases) this subsumes the above requirement, since $d\geqslant a_j\geqslant2$ automatically.
                        The reason for this requirement is simply that, if $a_j>d$ for some $j$, then we won't have any $x_j$ terms in $f$, and so we don't really have a `proper' polynomial in all the $x_j$.
                    }
                    \item if $a_i\nmid d$ then there exists some $j\neq i$ such that $a_i\mid(d-a_j)$.\footnote{%
                        This ensures that $f$ can satisfy the second condition: if $a_i\nmid d$ then it contains an $x_jx_i^{d-a_j}$ term.
                        In fact, this is arguably the most important condition -- we've already required that $f$ has no non-trivial non-constant factors (\cref{defn:plane-curves-in-wps}), and when we combine that with this requirement we see that we satisfy the hypotheses of \cite[8.4~Corollary]{IanoFletcher:2015wc}.
                        That is, we are asking that our plane curves be \emph{quasismooth}.
                    }
                \end{itemize}

                A plane curve $C=C_f\subset\pee$ is said to be \emph{sufficiently general} if its defining polynomial $f$ is sufficiently general.
            \end{definition}

            \emph{The main reason for \cref{defn:suff-general} is the \emph{quasismoothness} that it guarantees (see the footnotes), and that explains most of the inner workings of what follows from here on in.
            For a more thorough treatment, see \cite[Section~8]{IanoFletcher:2015wc}.}

            One of the reasons that we settle upon \cref{defn:suff-general} is because it gives us \cref{lem:when-is-p-i-in-c}, which will come in use later.
            But we also want to make sure that we haven't restricted ourselves so much that we end up studying a tiny subset of all possible plane curves.
            Another way of looking at \cref{defn:suff-general} is that, if we say a polynomial is of degree-$d$ then we want it to at least have an $x_i^{d/a_i}$ terms for all $i$.
            Since this might not always be possible, depending on the weighting, we sort of say that if this can't happen then we want the next best thing.
            \begin{lemma}\label{lem:when-is-p-i-in-c}
                Let $f\in\kathree$ be a sufficiently-general degree-$d$ weighted-homogeneous polynomial.
                Let $p_0=|1:0:0|$, $p_1=|0:1:0|$, and $p_2=|0:0:1|\in\pee$.
                Then $p_i\in C_f$ if and only if $a_i\nmid d$.
            \end{lemma}

            \begin{proof}
                By our definition of sufficiently general we have two cases:
                \begin{enumerate}[(i)]
                    \item if $a_i\mid d$ then there is an $x_i^{d/a_i}$ term in $f$, and so $f(p_i)\neq0$, and $p_i\not\in C_f$;
                    \item if $a_i\nmid d$ then every monomial containing a non-trivial power of $x_i$ also contains some non-trivial power of $x_j$ for $j\neq i$, and so every term of $f$ vanishes at $p_i$, thus $p_i\in C$.\qedhere
                \end{enumerate}
            \end{proof}

            Really though, condition (ii) in \cref{defn:suff-general} is superfluous if we assume that $f$ is also non-singular.
            The reason that we have that condition is to ensure that if $f$ is sufficiently general and if any of the $p_i$ are roots of $f$ then they are \emph{not} singular points.
            So if $p_i\in C$ then $C$ is not singular at $p_i$.
            But why do we make this restriction?

            The chain rule can tell us that $f$ is non-singular if $f(1,x,y)$, $f(x,1,y)$, and $f(x,y,1)$ are non-singular for all $(x,y)\in\aff^2\setminus\{(0,0)\}$ \emph{and} $f$ is non-singular at each of the $p_i$.
            So, by the above (which we state and prove in \cref{lem:p-i-not-sing-points}), it would suffice to check that $f$ is non-singular on the $U_i$ to show that $f$ is non-singular on the whole of $\pee$.
            This isn't a fact to which we appeal at all, and so we don't give all the gory details of using the chain rule, but it helps to reassure us slightly that our choice of definition might not be too bad -- if we have some non-singular affine curves that glue together to make a weighted projective curve then the resulting curve will also be non-singular.

            \begin{lemma}\label{lem:p-i-not-sing-points}
                Let $f\in\kathree$ be a sufficiently-general degree-$d$ weighted-homogeneous polynomial.
                If $p_i\in C$ (which happens if and only if $a_i\nmid d$) then $p_i$ is \emph{not} a singular point of $C$.
            \end{lemma}

            \begin{proof}
                Assume that $p_i\in C$, so that $a_i\nmid d$.
                By our definition then, $f$ contains an $x_jx_i^m$ term for some $j\neq i$, where $m=(d-a_j)/a_i$.
                So $\pdv{f}{x_j}$ contains an $x_i^m$ term.
                Further, every other monomial in $\pdv{f}{x_j}$ either has $x_j^l$ for some $l\geqslant1$ or is $x_k^l$ for some $l\geqslant1$.
                Either way, every other monomial vanishes at $p_i$ and the only remaining term is $x_i^m$, which evaluates to some non-zero scalar.
                Thus $\eval{\pdv{f}{x_j}}_{p_i}\neq0$, and so $p_i$ is not a singular point of $f$.
            \end{proof}

            Finally, we state and prove one more technical lemma here.
            This one is seemingly unrelated to anything we have mentioned so far, but ends up being a key part in the proof of \cref{thm:degree-genus}, so we get it out of the way now.
            The proof is messy but, in essence, simple.

            \begin{lemma}\label{lem:that-smoothness-lemma}
                Let $f\in\kathree$ be a sufficiently-general degree-$d$ non-singular weighted-homogeneous polynomial.
                Then $\cover{f}=\pi_\#(f)$ (as defined in \cref{lem:pi-hash-maps-homo-to-homo}) is such that all of
                \[
                    \cover{f}(0,1,\lambda),\quad\cover{f}(\lambda,0,1),\quad\cover{f}(1,\lambda,0)
                \]
                have no repeated roots when considered as polynomials in $\lambda$.

                {\color{purple}\textbf{Recall the additional hypothesis stated in \cref{lem:cover-is-smooth-pc-and-things}, which supposes that $\cover{f}$ is non-singular. This is not automatically true, given the hypotheses in this current lemma, and so must still be an \emph{additional} assumption.}}
            \end{lemma}

            \begin{proof}
                \emph{We prove that $\cover{f}(0,1,\lambda)$ has no repeated roots, and the other two claims follow in exactly the same manner, mutatis mutandis.}

                \bigskip

                Say for a contradiction that $\cover{f}(0,1,\lambda)$ has a multiple root $\lambda=c$, so that $(\lambda-c)^2\mid\cover{f}(0,1,\lambda)$.
                Also, since $\cover{f}(0,1,c)=0$ we know that $[0:1:c]\in C$.
                We aim to show that $[0:1:c]$ is a singular point of $\cover{f}$, contradicting the fact that $\cover{f}$ is non-singular.
                For easier reading we split this proof up into three parts -- one for each partial derivative of $\cover{f}$.

                \begin{enumerate}[(i)]
                    \item Now, if $y\neq0$ then $y^d\cover{f}(0,1,z/y)=\cover{f}(0,y,z)$, and so
                        \begin{align*}
                            \cover{f}(x,y,z) &= \cover{f}(0,y,z) + \sum_{i=1}^d c_i x^i y^j z^k \\
                            & = y^d\cover{f}(0,1,\lambda/y) + \sum_{i=1}^d c_i x^i y^j z^k \\
                            &= y^d\left(\frac{z}{y}-c\right)^2i(x,y,z) + xh(x,y,z) = (z-cy)^2g(x,y,z)+xh(x,y,z).
                        \end{align*}
                        where $c_i\in k\noz$ and $g,h,i$ are polynomials in $x,y,z$.
                        Thus
                        \[
                            \pdv{\cover{f}}{z} = (z-cy)g+(z-cy)^2\pdv{g}{z} + x\pdv{h}{z}
                        \]
                        and this evaluates to $0$ at $[x:y:z]=[0:1:c]$.
                    \item Next we examine two separate cases:
                        \begin{description}
                            \item[if $c\neq0$ then] $(z-cy)^2=c^2(y-z/c)^2$, and so $(y-z/c)\mid\pdv{\cover{f}}{y}$;
                            \item[if $c=0$ then] $\cover{f}=z^2g+xh$, and so $\pdv{\cover{f}}{y}=z^2\pdv{g}{y}+x\pdv{h}{y}$.
                        \end{description}
                        In both cases we see that $\pdv{\cover{f}}{y}$ evaluates to $0$ at $[x:y:z]=[0:1:c]$.
                    \item Finally,
                        \[
                            \pdv{\cover{f}}{x} = (z-cy)^2\pdv{g}{x} + x\pdv{h}{x} + h.
                        \]
                        But $\cover{f}(0,1,c)=0$ tells us that
                        \[
                            0=\eval{\bigg[\underbrace{(z-cy)^2g(x,y,z)}_{=0\text{ at }[x:y:z]=[0:1:c]}+h(x,y,z)\bigg]}_{[x:y:z]=[0:1:c]}
                        \]
                        so $h(0,1,c)=0$.
                        Thus $\pdv{\cover{f}}{x}$ evaluates to $0$ at $[x:y:z]=[0:1:c]$.\qedhere
                \end{enumerate}
            \end{proof}


        \subsubsection{Riemann-Hurwitz and the usual degree-genus formula} 
        \label{ssub:riemann_hurwitz_and_the_usual_degree_genus_formula}

            Here is where all the hard work in \cref{sub:weighted_projective_plane_curves_as_riemann_surfaces} pays of.
            With our idea of constructing a Riemann surface $\cover{C}\subset\pee^2$ for a plane curve $C\subset\pee$, along with a surjective map of Riemann surfaces $\pi\colon\cover{C}\twoheadrightarrow C$ with particularly nice properties (mainly the fact that it is holomorphic, but all of \cref{cor:quotient-set-up-with-pi-etc} comes in useful), we can now look at using one of the particularly powerful theorems from the study of Riemann surfaces: the \emph{Riemann-Hurwitz formula}.

            \begin{theorem}[Riemann-Hurwitz formula]\label{thm:rh-formula}
                Let $R,S$ be compact Riemann surfaces and $f\colon R\to S$ be a non-constant holomorphic map.
                Then
                \[
                    2g_R-2=\deg f(2g_S-2)-b(f)
                \]
                where $b(f)$ is the branching index
                \[
                    b(f) = \sum_{s\in S}\big(\deg f - |f^{-1}(s)|\big) = \sum_{s\in S}\left(\sum_{r\in f^{-1}(s)}(v_f(r)-1)\right)
                \]
                and $g_X=\frac12(\chi(X)+2)$ is the genus of the Riemann surface $X$.
            \end{theorem}

            \begin{proof}
                See almost any text on Riemann surfaces, e.g. \cite[Chapter~II,~Theorem~4.16]{Miranda:1995uz}.
            \end{proof}

            \begin{theorem}[Degree-genus formula for straight projective plane curves]\label{thm:straight-dg}
                Let $C\subset\pee^2$ be a non-singular degree-$d$ plane curve.
                Then $C$ is a Riemann surface with genus $g_C$ given by
                \[
                    g_C = \frac{(d-1)(d-2)}{2}.
                \]
            \end{theorem}

            \begin{proof}
                The fact that $C$ is a Riemann surface has already been proved in \cref{lem:str-proj-plane-curve-rs}, and \cite[Corollary~4.19]{Kirwan:1992wj} gives us the degree-genus formula.
            \end{proof}

            So we take our non-singular plane curve $C\subset\pee$ given by some sufficiently-general degree-$d$ weighted-homogeneous polynomial $f\in\kathree$, construct its straight cover $\cover{C}\subset\pee^2$ along with a quotient map $\pi\colon\cover{C}\twoheadrightarrow C$.
            Now $\cover{C}$ is also a non-singular plane curve defined by a homogeneous polynomial of degree-$d$, so we have our usual degree-genus formula for $\cover{C}$.
            But then the Riemann-Hurwitz formula tells us the genus of $C$ in terms of its degree $d$.
            If we write this all out properly then we expect to get some sort of degree-genus formula for non-singular sufficiently-general plane curves in $\pee$, and that is exactly what we get.

            \begin{theorem}[Degree-genus formula]\label{thm:degree-genus}
                \emph{Compare and contrast with \cite[Theorem~12.2]{IanoFletcher:2015wc}.}
                Let $C=C_f\subset\pathree$ be a non-singular plane curve where $f$ is weighted-homogeneous of degree $d$ and sufficiently general, in the sense of \cref{defn:suff-general}.
                Assume further that the straight cover $\cover{C}$ is non-singular (cf. \cref{lem:cover-is-smooth-pc-and-things}).
                Then, using the map $\pi$ as defined in \cref{cor:quotient-set-up-with-pi-etc},
                \[
                    g_C = \frac{1}{a_0a_1a_2}\left(\frac{(d-1)(d-2)}{2}-\left[\frac{b(\pi)}{2}+1-a_0a_1a_2\right]\right)
                \]
                where the branching index $b(\pi)$ is given by
                \[
                    b(\pi) = (d-1)\sum_{i=1}^3 (a_i-1) + \sum_{i=1}^3
                    \begin{cases}
                        a_i-1 &\text{ if }a_i\mid d;\\
                        a_0a_1a_2-1&\text{ if }a_i\nmid d.
                    \end{cases}
                \]
            \end{theorem}

            \begin{proof}
                First we appeal to \cref{cor:quotient-set-up-with-pi-etc}.
                This, along with the Riemann-Hurwitz formula and the degree-genus formula for straight projective plane curves (\cref{thm:rh-formula,thm:straight-dg}), tells us that
                \begin{equation}\label{eqn:thing-in-proof-to-rearrange}
                    2\frac{(d-1)(d-2)}{2}-2 = a_0a_1a_2(2g_C-2) + b(\pi).                
                \end{equation}

                Now \cite[Chapter~III,~Corollary~3.6]{Miranda:1995uz} tell us that
                \[
                    b(\pi) = \sum_{p\in\cover{C}}(|G_p|-1) = \sum_{i=1}^k\left(\frac{a_0a_1a_2}{|\pi^{-1}(y_i)|}-1\right)
                \]
                where $y_1,\ldots,y_k$ are \emph{all} the branch points of $\pi$.
                So all that remains to do is find and classify all of the branch points of $\pi$.

                Where would we expect to find branch points?
                Well, after a little bit of thinking, we see that the only branch points are those who have some zero coordinate, since that is the only way that the size of the preimage can drop.
                That is,
                \[
                    \bigg|\underbrace{\big\{[\sigma_0 x_0^{1/a_0}:\sigma_1 x_1^{1/a_1}:\sigma_2 x_2^{1/a_2}] \mid \sigma_i\in\mu^{a_i}\big\}}_{\pi^{-1}(|x_0:x_1:x_2|)}\bigg| < \underbrace{a_0a_1a_2}_{\deg\pi} \iff x_i=0\text{ for some }i
                \]
                where we once again use the assumption that $a$ is well-formed.\footnote{%
                    Or we could prove this in the other direction: since our field is algebraically closed and all the $a_i$ are pairwise coprime we know that $|\pi^{-1}(p)|=a_0a_1a_2$ for a point $p=|p_0:p_1:p_2|$ with $p_0,p_1,p_2\neq0$.
                }
                The points with some zero coordinates (i.e. the branch points $y_i$) split into two disjoint types: those with just one zero coordinate, and those with two.
                We introduce some temporary notation:\footnote{%
                    Here we use the convention that $x_3=x_0$.
                }
                \[
                    G_i = \{|x_0:x_1:x_2|\in C : x_i=0,~x_{i+1},x_{i+2}\neq0\}.
                \]
                So with $p_i$ as defined in \cref{lem:when-is-p-i-in-c} we can partition all of the branch points $y_i$ into four disjoint sets
                \[
                    \{y_1,\ldots,y_k\} = G_0\sqcup G_1\sqcup G_2\sqcup\{p_0,p_1,p_2\}.
                \]

                \textbf{We look first at the $G_i$:}
                we can see that if $x\in G_i$ then $|\pi^{-1}(x)|=a_{i+1}a_{i+2}=a_0a_1a_2/a_i$.
                So we know how each of the points in $G_i$ contribute to $b(\pi)$ and we are only left wondering \emph{how many} points there are in each $G_i$.
                Let's consider the example of $G_0$.
                Without loss of generality we can write points in $G_0$ as $|0:1:\lambda|$.
                Then asking how many points there are in $G_0$ is equivalent to asking how many \emph{non-zero} roots the polynomial $g_0(\lambda)=f(0,1,\lambda)$ has.
                The fundamental theorem of algebra tells us that it has $d$ roots overall, but counting multiplicity.
                However, \cref{lem:that-smoothness-lemma} tells us that all of the roots are distinct, and thus $g_0$ has $d$ distinct roots.
                By definition, $\lambda=0$ is a root if and only if $p_i\in C$ if and only if $a_i\nmid d$, and so we see that
                \[
                    |G_i| =
                    \begin{cases}
                        d &\text{if }a_i\mid d; \\
                        d-1 &\text{if }a_i\nmid d.
                    \end{cases}
                \]

                \textbf{Then we look at the $p_i$:}
                we don't need to worry about counting how many $p_i$ there are, since \cref{lem:when-is-p-i-in-c} tells us that $p_i\in C$ if and only if $a_i\nmid d$, and we can see that\footnote{%
                    Since, for example, $[\sigma:0:0]=[\rho:0:0]$ for any $\sigma,\rho\in\mu^{a_0}$.
                } $|\pi^{-1}(p_i)|=1$.

                Putting this all together we see that
                \begin{align*}
                    b(\pi) &= (d-1)\sum_{i=1}^3\left(\frac{a_0a_1a_2}{a_{i+1}a_{i+2}}-1\right) + \sum_{i=1}^3
                    \begin{cases}
                        \frac{a_0a_1a_2}{a_{i+1}a_{i+2}}-1 &\text{ if }a_i\mid d;\\
                        \frac{a_0a_1a_2}{1}-1&\text{ if }a_i\nmid d.
                    \end{cases} \\
                    &=
                    (d-1)\sum_{i=1}^3(a_i-1) + \sum_{i=1}^3
                    \begin{cases}
                        a_i-1 &\text{ if }a_i\mid d;\\
                        a_0a_1a_2-1&\text{ if }a_i\nmid d.
                    \end{cases}
                \end{align*}
                as claimed.
                All that then remains to prove the formula is to rearrange \cref{eqn:thing-in-proof-to-rearrange} into the given form.
            \end{proof}

            An interesting side-effect of \cref{thm:degree-genus} is that the complicated looking formula must always give an integer whenever $d$ and the $a_i$ satisfy the hypotheses, because we know that the genus of a Riemann surface is going to be an integer.
            This is similar to how the easiest way to show that $\frac{n!}{(n-k)!k!}$ is an integer is to note that $\binom{n}{k}$ is a way of counting things, and so \emph{must} be an integer.
        



\section{The view from where we've ended up} 
\label{sec:the_view_from_where_we_ve_got_to}

    \emph{This final section is disappointingly\footnote{%
            Unless you haven't been enjoying this paper, in which case this might be a welcome fact.
        } brief and also not entirely rigorous in the sense that, quite often we simplify things and hope that the reader will consult any of the given references before believing too much in anything read here.
        Also, some terminology or notation might not be explained.
        If this is ever the case then it is because it is considered standard (or as standard as mathematical notation can ever be) in the general literature of the subject, so any of the references provided, or other `classics', should clear up what it means.}

    \bigskip

    As my supervisor once said to me, `no piece of work is ever complete', and that is particularly true here.
    We have only scraped the surface of weighted projective space and its varieties, and we have done so in often simple language, which might not be the most natural way of explaining things (it often isn't).
    But there are a few things that came up during the writing of this text that the author found exciting and hopes that you might too.
    We present them to you in this section in an attempt to entice and lure people in to this exciting field of mathematics that comes under the vast umbrella that is `algebraic geometry'.

    During this text we have climbed a hill which, although minuscule in comparison to the towering peaks and ranges of mathematics as a whole, is not a hill to be sniffed at.
    From the top of this hill we can see a little bit more of the surrounding maths than before.
    It's still looming above us, but is ever so slightly more in focus and seems just that little bit more tangible and achievable.
    Let's take a look at the view from up here.

    \subsection{Elliptic curves and friends} 
    \label{sub:elliptic_curves_and_friends}

        Let $D$ be an ample divisor\footnote{%
            For our purposes we will only really look at point divisors, that is, $D=kp$ for some specific point $p\in C$ and $k\in\mathbb{N}$, and these are all ample, and in fact very ample for $k\geqslant2g_C+1$ by \cite[Chapter~IV,~Corollary~3.2(b)]{Hartshorne:1977we}.
        } on some plane curve $C\subset\pee$.
        Then we can define a graded ring
        \[
            R(C,D) = \bigoplus_{n\geqslant0}\mathcal{L}(C,nD)
        \]
        where $\mathcal{L}(C,nD)$ (written as just $\mathcal{L}(nD)$ when it is clear that we are working on $C$) is the Riemann-Roch space of meromorphic functions on $C$ with poles no worse than $nD$.
        Using reasonably standard notation we define $\mathcal{L}(C,nD)$ by:
        \[
            \mathcal{L}(C,nD) = \{f\colon C\to\cc \mid f\text{ is meromorphic and }(f)+D\geqslant0\}.
        \]

        In a loose sense, what it means for a divisor to be ample is that we can reconstruct the curve $C$ from this graded ring $R(C,D)$.
        To be slightly more precise,
        \[
            \proj R(C,D)\cong C,
        \]
        but the embedding of $C$ given by this $\proj$ construction might not naturally sit inside $\pee$.
        This is why it becomes very useful to look at this idea, since we get different ways of representing the same curve in different weighted projective spaces.

        \begin{example}[Elliptic curves]\label{ex:elliptic-first}
            \emph{There is an exercise sheet \cite{Reid:5t8DB_-1} by Miles Reid on graded rings which covers all of the below, and far far more besides.
            It should be reasonably easy to find online.}

            Let $C\subset\pee^2$ be given by a non-singular homogeneous cubic $f\in k[x_0,x_1,x_2]$, and define the divisor $D=p$ for some point $p\in C$.
            We can now use the Riemann-Roch theorem which tells us the dimension $\ell(nD)$ of $\mathcal{L}(nD)$ with a correction term involving a canonical divisor $\kappa$ of $C$:
            \[
                \ell(nD) - \ell(\kappa-nD) = \deg D + 1 - g.
            \]

            We know that $\deg nD = n$, since $D$ is just the point $p$, and we also know the genus $g=1$.
            It is a useful fact that $\deg\kappa=2g-2=0$, so here if $n\geqslant1$ then $\deg(\kappa-nD)<0$, and thus the correction term $\ell(\kappa-nD)=0$ disappears.
            Finally, we know that the only holomorphic functions on $C$ are the constant ones, and thus $\mathcal{L}(0)=\cc$.
            So Riemann-Roch tells us that
            \begin{equation}\label{eqn:riemann-roch-example}
                \ell(np) =
                \begin{cases}
                    1 & n=0\\
                    n & n\geqslant1.
                \end{cases}
            \end{equation}
            Using this fact, we can try to construct $R(C,np)$ a little bit more explicitly.

            \bigskip

            We've already said that $\mathcal{L}(0)=\cc$, so let's look at when $n\geqslant1$, using \cref{eqn:riemann-roch-example} to tell us how many elements we need for a basis\footnote{%
                Here we skim over issues of linear independence, in keeping with our theme this section of being quick and not very rigorous.
            }:
            \begin{description}
                \item[$n=1$:] Let $x\in\mathcal{L}(p)$ be such that $\langle x\rangle=\mathcal{L}(p)$.
                    Since $\ell(p)=1$ we know that $\mathcal{L}(p)\cong\cc$ and so we can take $x$ to be the image of $1\in\cc$ under this isomorphism.
                    Thus $x$ is an identity map;
                \item[$n=2$:] By definition, $x^2\in\mathcal{L}(2p)$, since
                    \[
                        (x^2)+2p = (x)+(x)+p+p = (x)+p+(x)+p \geqslant0.
                    \]
                    But $x$ is the identity, and so $x^2$ is really just a copy of $x$ that naturally lives inside $\mathcal{L}(2p)$, and hence still the identity.
                    So let $y$ be such that $\langle x,y\rangle=\mathcal{L}(2p)$;
                \item[$n=3$:] As above, $x^3,xy\in\mathcal{L}(3p)$, but we need one more basis element, so let $z$ be such that $\langle x,y,z\rangle=\mathcal{L}(3p)$;
                \item[$n=4$:] Now $x^4,x^2y,y^2,xz\in\mathcal{L}(4p)$, and no other combinations of $x,y,z$ are, so we have exactly the right amount of elements for a basis;
                \item[$n=5$:] Here $x^5,x^3y,xy^2,x^2z,yz\in\mathcal{L}(5p)$ gives us exactly the right amount of elements again;
                \item[$n=6$:] Things start to change in this case: $x^6,x^4y,x^2y^2,y^3,x^3z,xyz,z^2\in\mathcal{L}(6p)$, which is one element too many, so we must have some linear dependence between them all.
                    Since $x^k$ is an identity map, we see that the only `new' elements that we've got in the $n=6$ case are $y^3$ and $z^2$ (the only ones without a non-trivial power of $x$), and so the linear dependence must involve them.
                    By doing some linear change of coordinates\footnote{%
                        Full details in \cite[Chapter~IV,~Proposition~4.6]{Hartshorne:1977we}, or use $z\mapsto z+\alpha_3(x,y)$ and $y\mapsto y+\beta_2(x)$ with suitable constants, where the subscript represents the degree of the polynomial.
                    } we can assume that this equation is of the form
                    \begin{equation}\label{eqn:linear-relation-elliptic-curve}
                        z^2 = y^3+ax^4y+bx^6.
                    \end{equation}
                \item[$n\geqslant7$:] By doing some dimension counting we can show that $x,y,z$ generate $R(C,p)$ with only one relation between them, namely \cref{eqn:linear-relation-elliptic-curve}.
                    This sort of crosses over to using the Hilbert polynomial and some theory of generating functions, which is covered more in the next section, but the general idea is that 
                    \[
                        \big|\{x^\alpha y^\beta z^\gamma\in\mathcal{L}(np)\}\big| = \big|\{(\alpha,\beta,\gamma)\in(\nn\cup\{0\})^3 : \alpha+2\beta+3\gamma=n\}\big|
                    \]
                    which is exactly\footnote{%
                        See \cite[Section~3.15]{Wilf:2005to} for this specific problem, and the book as a whole for a great introduction to the theory and applications of generating functions.
                    } the $t^n$ coefficient of the series expansion of
                    \[
                        \frac{1}{(1-t)(1-t^2)(1-t^3)}.
                    \]
            \end{description}

            \bigskip

            So this tells us that
            \[
                R(C,D) \cong \frac{k_{(1,2,3)}[x,y,z]}{(g_6)}
            \]
            where $g_6=z^2-y^3+ax^4y+bx^6$ (from \cref{eqn:linear-relation-elliptic-curve}).
            Thus
            \[
                C \cong \proj R(C,D) \cong \van(g_6)\subset\pee(1,2,3)
            \]
            is an embedding of $C$ as a degree-$6$ plane curve in $\pee(1,2,3)$.

            \bigskip

            We can repeat this story but with $D=2p$, as this is still a very-ample divisor.
            But, up to a constant multiple of the grading (or not even that, depending on our grading convention for truncations), this is the same as looking at the $2$-nd truncation:
            \[
                R(C,2p) = \bigoplus_{n\geqslant0} \mathcal{L}(2np) = \bigoplus_{\substack{m\geqslant0\\2\mid m}} \mathcal{L}(mp) \cong \bigg(\bigoplus_{m\geqslant0}\mathcal{L}(mp)\bigg)^{(2)} = R(C,p)^{(2)}.
            \]
            This then gives us a different embedding, namely
            \[
                R(C,p)^{(2)} \cong \frac{k_{(1,1,2)}[x_0,x_1,y]}{(g_4)}
            \]
            which induces an embedding of $C$ as a degree-$4$ plane curve in $\pee(1,1,2)$.\footnote{%
                For details, \cite[Section~12.6]{IanoFletcher:2015wc} provides a thorough explanation.
            }

            But why stop at $2p$?
            It seems like we might as well look at $kp$ for $k\in\nn$.
            \cref{tab:elliptic-curves-truncations} summarises what we get when $D=kp$ for different values of $k$.\footnote{%
                Compare with the list in \cite[Section~12.4]{IanoFletcher:2015wc}.
            }
            As a reassuring fact, we see that \cref{thm:degree-genus} tells us that the genus is $1$ for all of these curves, as it should.

            \begin{table}[ht]
                \centering
                \begin{tabular}{lccl}
                    \toprule
                    $k$ & degree of curve(s) & ambient space & comments \\
                    \midrule
                    $1$ & $6$ & $\pee(1,2,3)$ & the classical Weierstrass equation \\
                    $2$ & $4$ & $\pee(1,1,2)$ & double cover of $\pee^1$ with $4$ branch points\\
                    $3$ & $3$ & $\pee(1,1,1)$ & plane cubic with an inflexion at infinity \\
                    $4$ & $2,2$ & $\pee(1,1,1,1)$ & intersection of two quadrics \\
                    $5$ & $5$ & $\pee(1,1,1,1,1)$ & $\pee^4$ section of $\operatorname{Grass}(2,5)\subset\pee^9$ \\
                    \bottomrule
                \end{tabular}
                \caption{Different embeddings of an elliptic curve $C$ coming from Veronese truncations of $R(C,p)$.}
                \label{tab:elliptic-curves-truncations}
            \end{table}

            We know\footnote{%
                \cite[Chapter~IV,~Corollary~3.2(b)]{Hartshorne:1977we} again.
            } that for $k\geqslant3$ the divisor $D_k$ is very ample, and we thus\footnote{%
                See the paragraph under the heading \emph{Ample Invertible Sheaves} in \cite[Chapter~II]{Hartshorne:1977we}, just before Remark~7.4.1.
            } get an embedding of $C$ into $\pee^{k-1}$, but as $k$ increases the associated description of $C$ gets more and more complicated.
            But we also know that a smooth\footnote{%
                For example, $C_7=\van(x^7+y^2z+xz^2)$.
            } $C_7\subset\pee(1,2,3)$ has genus $1$, and is thus also an elliptic curve.\footnote{%
                This is given as a brief example in \cite[Note~6.15(i)]{IanoFletcher:2015wc}.
            }
            Yet it doesn't appear in \cref{tab:elliptic-curves-truncations} at all, and we've just said that for $k\geqslant3$ all of the embeddings will be in $\pee^{k-1}$.
            So $C_7$ never arises from the above method.
        \end{example}

        It is an interesting question to ask \emph{why} $C_7$ never appears from the process in \cref{ex:elliptic-first}.
        Unfortunately, it is at this point that the author must once more throw up their hands in confession and admit uncertainty.
        A believe answer is that this specific embedding is not \emph{projectively normal}, but the details are beyond this text.


    \subsection{Syzygies and some homological algebra} 
    \label{sub:syzygies_and_some_homological_algebra}

        \emph{This section raises far more questions than it answers.
        Whether this is due to a lack of time or a lack of knowledge on the behalf of the author is left intentionally ambiguous.}

        \bigskip

        In \cref{ex:elliptic-first} we constructed a surjective graded ring homomorphism $\kathree\to R(C,D)$ with kernel $(g_6)$.
        Let's consider a more general situation where we have a surjective graded ring homomorphism $\vartheta\colon S_a=\kazn\to R(D)$ with kernel $(g_1,\ldots,g_k)$, where $\deg g_i=d_i$.
        Then we can write the sequence
        \begin{equation*}
            \begin{tikzcd}
                S_a \rar{\vartheta}
                & R(D) \arrow{r}
                & 0
            \end{tikzcd}
        \end{equation*}
        which is exact at $R(D)$, but not at $S_a$.
        To make this sequence exact we need a graded ring $S'$ along with a graded ring homomorphism $\vartheta'\colon S'\to S_a$ such that $\im\vartheta'=\ker\vartheta=(g_1,\ldots,g_k)$.
        Consider the map $\vartheta'\colon S'=\bigoplus_{i=1}^k S_a\to S_a$ given by
        \[
            \vartheta'\colon(f_1,\ldots,f_k)\mapsto f_1g_1+\ldots+f_kg_k
        \]
        which we write in matrix-like notation as $f\mapsto (g_1,\ldots,g_k)f$.
        This is not a graded ring homomorphism, but if we give $S'_a$ a different grading then we claim that we can make it one.
        Write $S_a[-d_i]$ to mean the graded ring $S_a$ but with a shift of grading by $-d_i$.
        That is, if $f\in S_a$ is such that $\deg f=d$ then when we consider $f$ as an element of $S_a[d_i]$ it has degree $d-d_i$.
        Let $S'=\bigoplus_{i=1}^k S_a[-d_i]$, so that we have the sequence
        \begin{equation*}
            \begin{tikzcd}
                \bigoplus_{i=1}^k S_a[-d_i] \rar{\vartheta'}
                & S_a \rar{\vartheta}
                & R(D) \arrow{r}
                & 0
            \end{tikzcd}
        \end{equation*}
        which is exact at $S_a$ and $R(D)$, but not at $\bigoplus_{i=1}^k S_a[-d_i]$.

        The question is, if we carry on finding $S'',S''',\ldots$ and $\vartheta'',\vartheta''',\ldots$ in a similar way, will we ever end up with an exact sequence of finite length?
        The answer is, maybe surprisingly, yes.

        \begin{theorem}[Hilbert syzygy theorem {(\cite[Theorem~1.13]{Eisenbud:1995tm})}]
            Let $R$ be a finitely-generated $k$ algebra.
            Then every finitely-generated $R$-module has a finite graded free resolution of length no more than $r$, by finitely-generated free modules.
        \end{theorem}

        The theorem is a bit stronger than what we've been asking for in our simplified language, but in essence it tells us that we can \emph{always} find a finite free resolution (i.e. turn our chain into an exact chain of finite length).
        Let's look at some simple examples.

        \begin{example}[$k=1$]\label{ex:k-equals-1}
            If we take $k=1$ then we can stop where we stopped above.
            That is, we have $\vartheta\colon S_a\to R(D)$ with $\ker(\vartheta)=(g)$ where $\deg g=d$.
            Then we have the (short) exact sequence
            \begin{equation*}
                \begin{tikzcd}
                    0 \arrow{r}
                    & S_a[-d] \rar{(g)}
                    & S_a \rar{\vartheta}
                    & R(D) \arrow{r}
                    & 0
                \end{tikzcd}\qedhere
            \end{equation*}
        \end{example}

        \begin{example}[$k=2$]\label{ex:k-equals-2}
            When $k=2$ we get a slightly more interesting case.
            Using the same matrix-like notation as before, we obtain the exact sequence
            \begin{equation*}
                \begin{tikzcd}
                    0 \arrow{r}
                    & S_a[-(d_1+d_2)] \rar{\smqty(g_2\\-g_1)}
                    & S_a[-d_1]\oplus S_a[-d_2] \rar{(g_1,g_2)}
                    & S_a \rar{"\vartheta"}
                    & R(D) \arrow{r}
                    & 0
                \end{tikzcd}
            \end{equation*}
            where $\smqty(g_2\\-g_1)\colon f\mapsto (g_2f, -g_1f)$.
        \end{example}

        This is all very interesting, but doesn't yet seem to relevant to what we were doing in the last section.
        But we now provide a few examples and some explanation as to why this is actually all very interesting.
        Compare and contrast the following with \cite[Proposition~4.3,~Example~4.4]{Reid:2002uy}.

        \begin{example}
            We see that a variety which falls under the category of that in \cref{ex:k-equals-1} (so given by $\van(g)\subset\pazn$) has the Hilbert series
            \[
                \frac{1-t^d}{\prod(1-t^{a_i})}
            \]
            and one in the same family as in \cref{ex:k-equals-2} has Hilbert series
            \[
                \frac{1-t^{d_1}-t^{d_2}+t^{d_1+d_2}}{\prod(1-t^{a_i})} = \frac{(1-t^{d_1})(1-t^{d_2})}{\prod(1-t^{a_i})}.
            \]
        \end{example}

        So if we can calculate the Hilbert series of a variety then, by writing it as a fraction and changing the denominator, we can get a rough idea of what embeddings we can get from it.
        We can see what the weights $a_i$ the denominator gives us, and we can say what sort of relations and syzygies its defining polynomials will satisfy by looking at the numerator.
        Why is this?

        Well, let $X\subset\pee^n$ be a projective variety with homogeneous coordinate ring $S(X)=k[x_0,\ldots,x_n]/\ide(X)$.
        Then we can construct the Hilbert series $P_X(t)$ of $X$ by considering some ample divisor $D$ on $X$ as we did above:
        \[
            P_{X,D}(t) = \sum_{m\geqslant0} \ell(mD)t^m.
        \]
        Using Riemann-Roch we can usually find $\ell(mD)$ explicitly, often by choosing $D=p$ to be a specifically nice point in $X$.
        Then, \emph{ignoring any issues of convergence} (by assuming that $|t|$ is small enough), we can try to rewrite this series as a single fraction.
        This is probably best explained here as an example, since we are being nowhere near rigorous enough to try to explain ourselves in proper mathematical language.

        \begin{example}[Elliptic curves, continued]\label{ex:elliptic-second}

            \emph{This example is a continuation \cref{ex:elliptic-first}, and so, in particular, all notation remains the same.}

            \bigskip

            We've already calculated $\ell(np)$ for $p\in C$, and so we know that the Hilbert polynomial of our elliptic curve is
            \[
                P_{C,p}(t) = \sum_{m\geqslant0}\ell(mp)t^m = 1+t+2t^2+3t^3+\ldots.
            \]
            Now, indeed
            \[
                \frac{1-t^6}{(1-t)(1-t^2)(1-t^3)} = 1+t+2t^2+\ldots = P_{C,p}(t)
            \]
            and by our previous comments this looks like an embedding of $C$ into $\pee(1,2,3)$ as the vanishing of a single degree-$6$ curve.
            Similarly,
            \begin{align*}
                P_{C,2p} &= 1+\ell(2p)t+\ell(4p)t^2+\ldots = \frac{1-t^4}{(1-t)(1-t)(1-t^2)} \\
                P_{C,3p} &= 1+\ell(3p)t+\ell(6p)t^2+\ldots = \frac{1-t^3}{(1-t)(1-t)(1-t)}
            \end{align*}
            which both agree with the embeddings that we already know. namely $C_4\subset\pee(1,1,2)$ and $C_3\subset\pee^2$.
        \end{example}

        The idea of generating functions and the combinatorics behind all this can really help to give us some intuition as to why we can read off such data, and why the numerator tells us about the relations.
        If $R(X,D)$ has generators $x_0,\ldots,x_n$ such that $x_i\in\mathcal{L}(X,a_iD)$ then to see how many combinations (i.e. products) of these lie in $\mathcal{L}(kD)$ we simply look at $c_k$, defined as the $k$-th coefficient in the series expansion of
        \[
            \frac{1}{\prod_i(1-t^{a_i})} = \sum_i c_i t^i.
        \]
        But if there are some relations between the $x_i$ then we won't be getting the full amount of distinct combinations all the time.
        That is, we won't have $c_k$ distinct combinations of the $x_i$ in $\mathcal{L}(X,c_kD)$ when $k$ becomes too large, since there will be some cancellation, and so we need to adjust our fraction to account for this -- we need to add in some higher order negative terms to reduce the $c_j$ for $j\geqslant k$.
        But then there might be some relations between the relations of the $x_i$, and in that case we have reduced the $c_j$ too much and we will actually have \emph{more} combinations of the $x_i$, and so we will need to put in some positive terms of higher order to make the $c_j$ larger, and so on.

        We end this paper with one final example, which is again more of an unfinished problem.

        \begin{example}
            Let's try to repeat the method that we applied in \cref{ex:elliptic-first} to study the genus $3$ Riemann surface given by a smooth quartic curve $C=C_4\subset\pee^2$, with the divisor $D=p\in C$.

            We know that a canonical divisor $\kappa$ on $C$ has degree $2g-2=4$ and is thus linearly equivalent to a hyperplane divisor $H=H_L$ for any line $L\subset\mathbb{P}^2$.
            So by Riemann-Roch, we know that $\ell(np)=n-2$ for $n\geqslant5$, since then $\deg(\kappa-nD)<0$ thus $\ell(\kappa-np)=0$.
            By definition, we can see that when $n=0,1$ we get $\ell(np)=1$,
            The question then is: how do we calculate $\ell(np)$ for $n=2,3,4$?
            We know that $\ell(np)$ is non-decreasing, and if we can further show that $\ell(np),\ell(\kappa-np)>0$ then we can use Clifford's theorem to obtain some bounds, but this still gives us a few possible options.

            It turns out that, actually, there is no one answer -- it depends on the point $p$ that we choose.
            We always have $\ell(2p)=1$ (by non-hyperellipticity), and we can actually choose $p$ to get any of the possible values of $\ell(3p)$ and $\ell(4p)$ that we like.\footnote{
                Have a look at a question asked by the author on \texttt{math.stackexhange.com}: \cite{Anonymous:5kIQ8xYm}.
            }
            So we have the following possibilities for the Hilbert series coefficients:
            \[
                \ell(np) =
                \begin{cases}
                    1,1,1,2,3,4,5,\ldots \\
                    1,1,2,2,3,4,5,\ldots \\
                    1,1,2,3,3,4,5,\ldots.
                \end{cases}
            \]
            For the sake of concreteness, let's just examine the first one here: assume that
            \[
                P_C(t) = 1+t+t^2+t^3+2t^4+3t^5+\ldots = 1+t+t^2+\sum_{m\geqslant4} (m-2)t^m.
            \]
            Following our naive approach from before, then, we find that we can pick elements
            \[
                v\in\mathcal{L}(p), w\in\mathcal{L}(4p), x\in\mathcal{L}(5p), y\in\mathcal{L}(6p), z\in\mathcal{L}(7p)
            \]
            that generate $\mathcal{L}(np)$ without any relations until we hit $n=10$, and then we have too many elements.

            The $n$-th coefficient of the series expansion of $\lambda(t)=1/(1-t)(1-t^4)(1-t^5)(1-t^6)(1-t^7)$ tell us how many elements of $\mathcal{L}(np)$ we can generate with $v,w,x,y,z$, as in \cref{ex:elliptic-first}.
            Since we have too many elements for $n\geqslant10$ we expect our Hilbert polynomial to have the same denominator as $\lambda$, but with a negative term of degree $10$ in the numerator to lower the coefficients for $t^n$ in $\lambda$ when $n\geqslant10$.
            This is where combinatoric intuition can help us understand a problem about projective plane curves.

            If we do some algebra, we see that
            \begin{align*}
                P_C(t) & = 1+t+t^2+t^3+2t^4+3t^5+\ldots \\
                & = \frac{1}{1-t} + t^4(1+2t+3t^2+\ldots) \\
                & = \frac{1}{1-t} + t^4\dv{t}(1+t+t^2+\ldots) \\
                & = \frac{1}{1-t} + t^4\dv{t}\left(\frac{1}{1-z}\right) \\
                & = \frac{1}{1-t} + \frac{t^4}{(1-t)^2} \\
                & = \frac{1-t+t^4}{(1-t)^2}
                = \ldots \\
                &= \frac{1-t^{10}-t^{11}-2t^{12}-t^{13}-t^{14}+t^{16}+2 t^{17}+2 t^{18}+2 t^{19}+t^{20}-t^{23}-t^{24}-t^{25}}{(1-t)(1-t^4)(1-t^5)(1-t^6)(1-t^7)}.
            \end{align*}
            From this we can see that our discovery that there must be a relation of degree $10$ seems correct, and we speculate (and \emph{only} speculate, not claim with any amount of certainty) by looking at how the signs change that there are also relations of degree $11,12,13,14$ with syzygies of degree $16,17,18,19,20$ and higher syzygies between them of degree $23,24,25$.
        \end{example}




\end{document}